\documentclass[reqno]{alea2}
\usepackage{natbib}
\usepackage{fancyhdr}
\usepackage{graphicx}

\usepackage{amssymb}
\usepackage{amsmath}

\renewcommand{\cite}{\citet}

\makeatletter \@addtoreset{equation}{section} \makeatother

\renewcommand\theequation{\thesection.\arabic{equation}}
\renewcommand\thefigure{\thesection.\@arabic\c@figure}
\renewcommand\thetable{\thesection.\@arabic\c@table}

\newtheorem{theorem}{Theorem}[section]
\newtheorem{lemma}[theorem]{Lemma}
\newtheorem{proposition}[theorem]{Proposition}

\newtheorem{definition}[theorem]{Definition}

\newtheorem{remarks}[theorem]{Remarks}
\newtheorem{example}[theorem]{Example}
\def\prf{\noindent{Proof~: }}

\newcounter{tictac1}
\newcounter{tictac2}
\newcounter{tictac3}
\newcounter{tictac4}
\newcounter{tictac6}
\newcounter{tictac7}
\newcounter{tictac8}
\newcounter{tictac9}
\newcounter{tictac10}

\newenvironment{fleuvea}{
   \begin{list}{\rm{$\textbf{(\alph{tictac2})}$} }{\usecounter{tictac2}
\leftmargin 1cm\labelwidth 2em}}{\end{list}}
\newenvironment{fleuveP}{
   \begin{list}{\rm{$\textbf{(P\arabic{tictac3})}$} }{\usecounter{tictac3}
\leftmargin 1cm\labelwidth 2em}}{\end{list}}

\newenvironment{fleuve11}{
   \begin{list}{\rm{$\textbf{(\arabic{tictac8})}$} }{\usecounter{tictac8}
\leftmargin 1cm\labelwidth 2em}}{\end{list}}
\newenvironment{fleuveII}{
   \begin{list}{\rm{$\textbf{(\Roman{tictac3})}$} }{\usecounter{tictac3}
\leftmargin 1cm\labelwidth 2em}}{\end{list}}
\newenvironment{fleuveIII}{
   \begin{list}{\addtocounter{tictac9}{3}\rm{$\textbf{(\Roman{tictac9})}$} }{
\leftmargin 1cm\labelwidth 2em}}{\end{list}}
\newenvironment{fleuveaa}{
   \begin{list}{\addtocounter{tictac10}{1}\rm{$\textbf{(\alph{tictac10})}$} }{
\leftmargin 1cm\labelwidth 2em}}{\end{list}}

\newcommand{\reff}[1]{(\ref{#1})}
\newcommand{\comp}[1]{#1^{\rm c}}

\newcommand{\zzz}{\mathbb{Z}}
\newcommand{\tribu}{\mathcal{F}}

\def\1{\rlap{\mbox{\small\rm 1}}\kern.15em 1}
\def\ind#1{\1_{#1}}
\def\build#1_#2^#3{\mathrel{\mathop{\kern 0pt#1}\limits_{#2}^{#3}}}
\def\tend#1#2#3{\build\hbox to 12mm{\rightarrowfill}_{#1\rightarrow
#2}^{#3}}

{}
\def\converge#1#2#3{\build\hbox to
15mm{\rightarrowfill}_{\hbox{\scriptsize #3}}^{#1\rightarrow #2}}
\def\converg#1#2#3{\build\hbox to
15mm{\rightarrowfill}_{\hbox{\scriptsize #3}}^{#1\uparrow #2}}
\def\embf#1{\emph{\bf #1}}

\begin{document}

\date{August 27, 2005, accepted November 9, 2006}
\keywords{specification, Gibbs-measure,
non-nullness, quasilocality.}
\subjclass{82B05.}

\author{Roberto Fern\'andez and Gr\'egory Maillard}
\address{Roberto Fern\'andez, Laboratoire de Math\'ematiques Rapha{\"e}l Salem,
UMR 6085 CNRS-Universit\'e de Rouen, Avenue de l'Universit\'e, BP 12,
F-76801 Saint \'Etienne du Rouvray, France}
\email{roberto.fernandez@univ-rouen.fr}
\urladdr{\tt http://www.univ-rouen.fr/LMRS/Persopage/Fernandez/index.html}

\address{Gr\'egory Maillard, Institut de Math\'ematiques,
\'Ecole Polytechnique F\'ed\'erale de Lausanne, Station 8, CH-1015 Lausanne,
Switzerland}
\email{gregory.maillard@epfl.ch}
\urladdr{\tt http://ima.epfl.ch/\~{}maillard/}

\title[Construction of a specification from its singleton part] {Construction of a specification from its singleton part}

\begin{abstract}
We state a construction theorem for specifications starting from single-site
conditional probabilities (singleton part). We consider general single-site
spaces and kernels that are absolutely continuous with respect to a chosen
product measure (free measure). Under a natural order-consistency assumption
and weak non-nullness requirements we show existence and uniqueness of the
specification extending the given singleton part. We determine conditions
granting the continuity of the specification. In addition, we show that,
within a class of measures with suitable support properties, consistency
with singletons implies consistency with the full specification.
\end{abstract}

\maketitle

\section{Introduction}

A \emph{specification} on a product space of the form
$E^{\zzz^{d}}$ is a family of probability kernels labelled by the
finite subsets $\Lambda\subset\zzz^{d}$ satisfying the
requirements of a consistent system of conditional probabilities.
They are the central objects of mathematical statistical mechanics,
see, for instance, Georgii (1988)\nocite{geo88}.  In this paper we
determine conditions that guarantee the (re)construction of a
specification from single-site conditional probabilities
(\emph{singletons}).  Such a scenario yields an interesting
simplification of the theory of specifications, and sets it in a
framework analogous to that of discrete-time stochastic processes,
traditionally defined and characterized by properties of single-site
transition probabilities.

The issue of singleton characterization of specifications stems
already from Dobrushin's (1968)\nocite{dob68} seminal work.  His
remarks were taken up by Flood and Sullivan
(1980)\nocite{flosul80} and lead to Theorem (1.33) in Georgii
(1988)\nocite{geo88}.  These references studied the reconstruction
problem, namely how to recover a pre-existing specification
starting from the singletons or from a subspecification.  More
recently, Dachian and Nahapetian (2001,
2004)\nocite{dacnah01,dacnah04}, and us (Fern\'andez and Maillard,
2004)\nocite{fermai03a} have addressed the more general
construction problem under complementary hypotheses. The key
issue, for these constructions, is the degree of ``nullness"
allowed for the specification, that is, the presence of ``excluded
configurations" leading to zero probability weights. The results
by Dachian and Nahapetian are suited to situations where the
exclusions come from asymptotic (measurable at infinity) events.
In contrast, our 2004 results apply for the case of local
exclusions (``grammars"). Our present results are an extension of
those by Dachian and Nahapetian, and coincide with our 2004
construction only for non-null singletons. Moreover, our proof,
while inspired in existing proofs, offers an alternative
formulation that, we believe, clarifies the algebraic and
measure-theoretical properties involved.

We work with general single-spin spaces and consider singletons
that are absolutely continuous with respect to a pre-established
product measure (\emph{free measure}).  This is the natural
framework from the physical point of view. We demand two key
conditions, besides the obvious finiteness and normalization
requirements: {\bf (H1)} some degree of non-nullness, and {\bf
(H2)} a compatibility condition.  The former is the extension, to
our framework, of Dachian and Nahapetian's (2004) \emph{very weak
positivity}.
Condition (H2) is the adaptation, in the absence of strict
positivity, of the compatibility identity (A.9) in Fern\'andez and
Maillard (2004). This is a partially integrated condition that,
for finite spins, coincide with the pointwise condition imposed by
Dachian and Nahapetian (2004) (defining what they call
\emph{1-point specifications}).  Under reasonable hypotheses, an
almost-sure version of (H2) follows from (H1) (Proposition
\ref{prop1} below).

Under these conditions we show (Theorem \ref{th1}) that there
exists a unique specification that is absolutely continuous with
respect to the free measure and whose single-site probabilities
coincide with the given singletons.  The proof provides a
recursive construction of this specification [formulas
\reff{r.1}--\reff{r.2}].  Our scheme makes no use of the possible
continuity of each singleton with respect to exterior
configurations.  As such, it is equally applicable to Gibbsian
(Kozlov, 1974\nocite{koz74} and Sullivan, 1973\nocite{sul73}) and
non-Gibbsian (van Enter, Fern\'andez and Sokal,
1993\nocite{vEFS_JSP}) theories.  Nevertheless we determine a
natural condition ensuring that the continuity of singletons lead
to a continuous specification.  Furthermore, in the third part of
our theorem, we establish a natural class of measures for which
consistency with the original singletons implies consistency with
the full specification constructed from them.  The validity of a
similar implication for general measures remains open in this
setting (it has been established for local exclusion rules in
Fern\'andez and Maillard, 2004).

We illustrate our results with a simple example showing the actual
meaning of the different hypotheses.   We also present a rather
detailed comparison of our theorem with the preceding results.

\section{Preliminaries}

We consider a general measurable space $(E,\mathcal{E})$ and the
product space $\Omega = E^{\zzz^{d}}$ for $d\ge 1$
(\emph{configuration space}), endowed with the product
$\sigma$-algebra $\tribu = \mathcal{E}^{\zzz^{d}}$.  Our notation will
be fairly standard.  We shall denote $\mathcal{P}(\Omega,\tribu)$ the
set of probability measures on $(\Omega,\tribu)$. Support sites will be
indicated with subscripts, if different from the whole of $\zzz^d$.
For example, if $U \subset \zzz^{d}$ we denote $\Omega_U = E^U$ and
$\tribu_U$ the sub-$\sigma$-algebra of $\tribu$ generated by the cylinders
with base in $\Omega_U$.  Likewise $\sigma\in\Omega$,
$\sigma_\Lambda\in\Omega_\Lambda$.  ``Concatenated'' configurations
will be denoted as customary: If $\Lambda, \Delta \subset \zzz^d$ are
disjoint, $x_\Lambda \sigma_\Delta\in\Omega_{\Lambda\cup\Delta}$ is
the configuration coinciding with $x_\Lambda$ on $\Lambda$ and with
$\sigma_\Delta$ on $\Delta$, while $x_\Lambda \sigma_\Delta\omega$ is
the configuration in $\Omega$ which in addition is equal to $\omega$
on $\comp{(\Lambda\cup\Delta)}$.  One-site sets will be labelled just
by the site, for instance we shall write $\omega_i$ instead of
$\omega_{\{i\}}$.  For each $U \subset \zzz^{d}$, its cardinal will be
denoted $\left|U \right|$, its indicator function $\ind{U}$ and the
set of its finite subsets $\mathcal{S}(U)$.  We shall abbreviate
$\mathcal{S}\triangleq \mathcal{S}(\zzz^{d})$.  For $\Lambda \subset
\zzz^d$ and $i\in \Lambda$ we denote $\Lambda_{i}^{*} \triangleq
\Lambda \setminus \{i\}$ ($|\Lambda| \geq 1$).
Throughout this paper we adopt the convention ``$1/\infty=0$''.

We recall that a \emph{measure kernel} on $\tribu\times \Omega$ is a
map $\gamma(\,\cdot\mid\cdot\,):\tribu\times \Omega \to \mathbb{R}$
such that $\gamma(\,\cdot\mid \omega)$ is a measure for each
$\omega\in\Omega$ while $\gamma(A\mid\cdot\,)$ is
$\tribu$--measurable for each event $A\in\tribu$.  If each
$\gamma(\,\cdot\mid \omega)$ is a probability measure the kernel is
called a \emph{probability kernel}.  To obtain cleaner formulas we
shall adopt operator-like notations to handle kernels.  Thus, for
kernels $\gamma$ and $\widetilde\gamma$ and non-negative measurable
functions $f$ and $\rho$, we shall denote:
\begin{itemize}
\item $\gamma(f)$ for the measurable function $\int
  f(\eta)\,\gamma(d\eta\mid\,\cdot\,)$.
\item $\gamma\, \widetilde\gamma$ for the composed kernel defined by
  $(\gamma\,\widetilde\gamma)(f) = \gamma(\widetilde\gamma(f))$.
\item $\rho\,\gamma$ for the product kernel defined by $(\rho\,
  \gamma)(f) =\gamma(\rho f)$.
\end{itemize}

The following are the only two definitions needed for this paper.

\begin{definition}\label{d.r1}
A \embf{specification} on $\left(\Omega,\tribu\right)$
is a family of probability kernels $\left\{ \gamma_{\Lambda}
\right\}_{\Lambda \in \mathcal{S}}$ such that for all $\Lambda$ in
$\mathcal{S}$,
\begin{fleuvea}
\item $\gamma_{\Lambda}(A \mid \cdot\, ) \in \tribu _{\comp{\Lambda}}$
for each $A \in \tribu$ .
\item $\gamma_{\Lambda}(B \mid \omega) = \ind{B}(\omega)$
for each $B \in \tribu _{\comp{\Lambda}}$ and
$\omega \in \Omega$.
\item For each $\Delta \in \mathcal{S}$ with $\Delta \supset
\Lambda$,
\begin{equation}
\label{gibbs5} \gamma_{\Delta}\,\gamma_{\Lambda} \;=\;
\gamma_{\Delta}\;.
\end{equation}
\end{fleuvea}
\end{definition}

The last property is called \emph{consistency}.  It is stronger than
the \emph{almost} consistency of the finite-volume conditional
probabilities of a measure on $(\Omega,\tribu)$.  Without further
requirements, this strengthening is usually illusory: If
$(E,\mathcal{E})$ is a standard Borel space, each measure on
$(\Omega,\tribu)$ is consistent with some specification (Sokal,
1981\nocite{sok81}).  Matters become more delicate if in addition
kernels are requested to be continuous with respect to their second
variable, that is if the Feller property is imposed. Consistency with
a continuous specification is the hallmark of Gibbsianness.  See, for
instance, van Enter, Maes and Shlosman (2000)\nocite{entmaeshl00} for
a survey of the different notions and issues arising when this
continuity is absent.

\begin{definition}\label{d.con}
A probability measure $\mu$ on $(\Omega, \mathcal{F})$ is said to be
\embf{consistent} with a specification
$\{\gamma_\Lambda\}_{\Lambda\in\mathcal{S}}$ if
\begin{equation}\label{gibbs9}
\mu \,\gamma_{\Lambda}\;=\;\mu \quad \text{for every }
\Lambda \in \mathcal{S}.
\end{equation}
The family of these measures will be denoted
$\mathcal{G}\bigl(\{\gamma_\Lambda\}_{\Lambda\in\mathcal{S}}\bigr)$.
\end{definition}

\section{Main hypotheses}\label{s.hyp}

Throughout this paper we fix a family $\left(\lambda^{i}\right)_{i \in
  \zzz^{d}}$ of \emph{a priori} (non-negative) measures on
$(E,\mathcal{E})$.  Its choice is in general canonically dictated by
the structure of the single-spin space $E$.  For instance, if $E$
admits some group structure all $\lambda^i$ are chosen equal to the
corresponding Haar measure (the ``more symmetric'' measure).  For each
$\Lambda\subset\zzz^d$, let $\lambda^{\Lambda} \triangleq
\bigotimes_{i \in \Lambda} \lambda^{i}$ (\emph{free measure on
  $\Lambda$}) and $\lambda_{\Lambda}$ denote the kernel (\emph{free
  kernel on $\Lambda$}) defined by
\begin{equation}\label{spe7}
\lambda_{\Lambda}(h \mid \omega) \;=\;
\left(\lambda^{\Lambda} \otimes \delta_{\omega_{\comp{\Lambda}}}
\right)(h) \;=\; \int h(\sigma_{\Lambda} \omega) \,
\lambda^{\Lambda}(d \sigma_{\Lambda})
\end{equation}
for every measurable function $h$ and configuration $\omega$.

Except in part {\bf (III)} of our Theorem \ref{th1}, the measures $\lambda^i$
are not required to be normalized or even finite.  The lack of normalization
is the only aspect that could prevent the family
$\left(\lambda_\Lambda\right)_{\Lambda\in\mathcal{S}}$ from being a
specification.  Indeed, this family satisfies (a) and (b) of
Definition \ref{d.r1} and, furthermore, the following factorization
property:
\begin{equation}\label{spe15}
\lambda_{\Lambda \cup \Delta} \;=\;
\lambda_{\Lambda} \lambda_{\Delta}\;,
\end{equation}
for each pair of disjoint sets $\Lambda,\Delta\subset\zzz^d$.  If the
kernels are normalized, this is a strengthening of the consistency
condition (c) above.

We shall construct specifications by multiplying each kernel
$\lambda_\Lambda$ by a suitable measurable function $\rho_\Lambda$.
The resulting kernels can be interpreted as \emph{dependent} or
\emph{interacting} kernels.  A family
$\left(\rho_\Lambda\right)_{\Lambda\in\mathcal{S}}$ yielding an
interacting kernel is called a $\lambda$-\emph{modification} in
Georgii's (1988)\nocite{geo88} treatise (see, specially, Section 1.3).
If $E$ is countable and each $\lambda^i$ is (a multiple of) the
counting measure, every specification is obtained in this form.

Our specifications will be built starting from a family of single-site
kernels of the form $\rho_i\lambda_i$, $i\in\zzz$.  The following
definitions state the crucial hypotheses granting the feasibility of
our construction.

\begin{definition}\label{d.h1}
  A family $\{\rho_i\}_{i \in \zzz^{d}}$, of $\tribu$-measurable
  functions $\rho_i:\Omega \to [0,\infty[$ satisfies \embf{hypothesis
  (H1)} if for each $\omega \in \Omega$, $j \in \zzz^{d}$ and
  $V \in \mathcal{S}\left(\comp{\{j\}}\right)$, there exists $x_{j} \in
  \Omega_{j}$ such that
\begin{equation}\label{spe16}
\rho_{j}\left(x_{j} \sigma_V \omega\right)>0, \quad \forall \,
\sigma_V \in \Omega_V\;,
\end{equation}
and, for every $i \in \zzz^{d} : i \neq j$,
\begin{equation}\label{spe17}
\inf \left\{\lambda_{i}\Bigr(\rho_{i} \, \rho_{j}^{-1}
  \Bigr)\left(x_{j} \sigma_{V} \omega \right) :
\sigma_{V} \in \Omega_{V} \right\}>0
\end{equation}
and
\begin{equation}\label{spe18} \sup
\left\{\lambda_{i}\Bigr( \rho_{i} \,
\rho_{j}^{-1} \Bigr)\left( x_{j} \sigma_{V} \omega \right) :
\sigma_{V} \in \Omega_{V} \right\}<\infty\;.
\end{equation}
We denote
\begin{equation}\label{spe19}
b(j ,V,\omega) \triangleq
\Bigl\{x_{j} \in \Omega_j \text{ satisfying
    \reff{spe16}--\reff{spe18}}\Bigr\}\;
\end{equation}
and
\begin{equation}\label{spe19*}
B(j,V) \; \triangleq\; \Bigl\{\omega\in\Omega : \omega_j\in
b(j,V,\omega)\Bigr\}\;.
\end{equation}
Furthermore for every $W\in\mathcal{S}(\comp{V})$,
\begin{equation}\label{spe13}
b(V, W, \omega) \; \triangleq\; \Bigl\{x_{V} \in
  \Omega_{V} : x_{k} \in b(k,V_{k}^{*} \cup W, \omega)
  \text{ for every } k \in V \Bigr\}\;.
\end{equation}
\end{definition}

If $E$ is finite, hypothesis (H1) is exactly the condition of very
weak positivity introduced by Dachian and Nahapetian (2004).  Our sets
$b(j, V, \omega)$ correspond to \emph{good} (``\emph{bonnes}'') configurations
at the site $j$ given $\omega$ outside $V\cup\{j\}$, while the $B(j,V)$
correspond
to configurations that are good in a more global sense. In both cases,
this ``goodness" must be uniform with respect to the configurations in $V$.
The product
structure of the sets $b(\Lambda, W, \omega)$, embodied in definition
\reff{spe13}, is
essential for our procedure and prevents its immediate extension to
other than product spaces.

\begin{definition}\label{d.h2}
  A family $\{\rho_i\}_{i \in \zzz^{d}}$, of $\tribu$-measurable
  functions $\rho_i:\Omega \to [0,\infty[$ satisfies \embf{hypothesis
(H2)} if for each $i,j$ in $\zzz^{d}$ and $\omega \in \Omega$, the
  following is true:

For each $x_i \in b(i, \{j\}, \omega)$ and $x_j \in b(j,
  \{i\}, \omega)$,
\begin{equation}\label{r.20}
\frac{\rho_{i}(\omega) \, \rho_j (x_i \omega)}{\rho_i (x_i \omega)\,
\lambda_{j}\left(\rho_{j} \, \rho_{i}^{-1}\right)(x_i \omega)} \;=\;
  \frac{\rho_{j}(\omega) \, \rho_i (x_j
\omega)}{\rho_j (x_j \omega)\, \lambda_{i}\left(\rho_{i} \,
  \rho_{j}^{-1}\right)(x_j \omega)} \;.
\end{equation}
\end{definition}
As a consequence, the map $R_{i}^{j}:\Omega \longrightarrow ]0,+\infty]$
defined
  by
\begin{equation}\label{spe20}
R_{i}^{j}(\omega) \;=\;
\left(\frac{\rho_{i}}{\rho_{j}} \times \lambda_{j} \left(\rho_{j} \,
  \rho_{i}^{-1}\right) \right) \left(x_{i} \omega \right)
\end{equation}
is independent of the choice of $x_{i}  \in b(i,\{j\},\omega)$ and
  hence defines a $\tribu_{\comp{\{i\}}}$-measurable map.

Let us pause to discuss the meaning and motivation of these
hypotheses.  The conditions \reff{spe16}--\reff{spe18} in (H1)
imply that the denominators in \reff{r.20} and the numerator in
\reff{spe20} are neither zero nor infinity.  The denominator can be
zero in the latter, in which case $R_i^j(\omega)=\infty$.

As the reader will see, $R_i^j$ is what is needed to fulfill the
identity
\begin{equation}\label{r.31}
\rho_{\{i,j\}}(\omega) \;=\; \frac{\rho_j(\omega)}{R_i^j(\omega)}\;.
\end{equation}
Due to the $i\leftrightarrow j$ symmetry of the LHS, this identity
must be accompanied by the consistency requirement

\begin{equation}\label{c.2}
\frac{\rho_i(\omega)}{R_i^j(\omega)}\;=\;
\frac{\rho_j(\omega)}{R_j^i(\omega)}\;.
\end{equation}
Under strict positivity hypotheses, identity \reff{r.31} holds with
\begin{equation}\label{r.30}
R_i^j(\omega) \;=\; \lambda_i \left(\rho_i \,
  \rho_j^{-1}\right) (\omega)\;,
\end{equation}
as exploited in Georgii (1988), Theorem (1.33), or in Fern\'andez and
Maillard (2004), Appendix.  The consistency condition \reff{c.2} is
imposed as a further hypothesis in the latter reference, while it is
automatic in the former because the singletons are known to come from
a specification in the first place.  A look to our arguments in the
aforementioned appendix convinced us that to extend them to weakly
positive cases we should at least start from the following desideratum:
\begin{itemize}
  \item[(i)] Identities \reff{r.31} and \reff{c.2} must be true.
  \item[(ii)] Definition \reff{r.30} must be verified whenever the RHS
    is meaningful.
  \item[(iii)] $R_i^j$ must be $\tribu_{\comp{\{i\}}}$-measurable [as
    in \reff{r.30}].
\end{itemize}
The quantity $\lambda_j\left(\rho_j \, \rho_i^{-1}\right) (\omega)$
is well defined whenever $\omega_i=x_i  \in b(i,\{j\},\omega)$.  In this
case the validity of \reff{c.2} and (iii) of the desideratum implies
\begin{equation}
\frac{\rho_i(x_i\omega)}{R_i^j(\omega)}\;=\;
\frac{\rho_j(x_i\omega)}{\lambda_j\left(\rho_j \,
  \rho_i^{-1}\right) (x_i\omega)}\;.
\end{equation}
This explains \reff{spe20}.  With this definition of $R_i^j$, identity
\reff{r.20} is exactly \reff{c.2}.  Our theorem below shows
that, in fact, the above desideratum is basically all that is needed to
make a successful construction.

\section{Results}\label{s.main}

\begin{theorem}\label{th1}
  Let $\{\rho_i\}_{i \in \zzz^{d}}$ be a family of $\tribu$-measurable
  functions $\rho_i:\Omega \to [0,\infty[$ satisfying:
\begin{fleuvea}
\item For every $i$ in $\zzz^{d}$,
\begin{equation}\label{spe33}
 \lambda_{i} (\rho_{i} \mid\omega)\;=\;1\;,
\end{equation}
for all $\omega \in \Omega$.
\item Hypotheses \rm{(H1)} and \rm{(H2)}.
\end{fleuvea}
Then there exists a family $\left\{ \rho_{\Lambda} \right\}_{\Lambda
  \in \mathcal{S}}$ of measurable functions $\rho_\Lambda:\Omega\to
[0,\infty[$, with $\rho_{\{i\}}=\rho_i$, such that the family of
kernels $\left\{ \rho_{\Lambda} \lambda_{\Lambda}\right\}_{\Lambda \in
  \mathcal{S}}$ is a specification.  Furthermore:
\begin{fleuveII}
\item If
\begin{equation}\label{r.34}
\lambda^j\bigl( b(j,V,\omega)\bigr)>0
\end{equation}
for each $\omega \in \Omega$, $V \in \mathcal{S}$ and $j \in
\comp{V}$, then there exists exactly one family $\left\{ \rho_{\Lambda}
\right\}_{\Lambda \in \mathcal{S}}$ with the above property.
\item Suppose that $E$ is a topological space and $\mathcal{E}$ its
  borelian $\sigma$-algebra, and consider the product topology for
  $\Omega$.  If the functions $\rho_i$ are sequentially continuous and
  for each $V \in \mathcal{S}$ and $j\in\comp{V}$ there exists
  $x_j \in \bigcap_{\omega} b(j,V,\omega)$
  such that
  \begin{equation}\label{cont1}
   \int \sup_{\omega}
   \Bigl[(\rho_i\,\rho_j^{-1})(\sigma_i x_j \omega)\Bigr]
    \,\lambda^i(d\sigma_i) \;<\;\infty
\end{equation}
for all $i \in V$, then the functions $\rho_\Lambda$,
 $\Lambda\in\mathcal{S}$, are sequentially continuous.
\end{fleuveII}

\noindent
Explicitly, the functions $\rho_\Lambda$ are recursively defined
throughout the identity
\begin{equation}\label{r.1}
\rho_{\Theta \cup \Gamma}(\omega) \;=\;
\frac{\rho_{\Theta}(\omega)}{R_{\Theta}^{\Gamma}(\omega)}\;,
\end{equation}
valid for every $\Theta \in \mathcal{S}$, $\Gamma \in
\mathcal{S}(\comp{\Theta})$, $\omega \in \Omega$,
where
\begin{equation}\label{r.2}
R_{\Theta}^{\Gamma}(\omega) \;=\;
\left(\frac{\rho_\Theta}{\rho_\Gamma} \times
  \lambda_{\Gamma}\left(\rho_{\Gamma} \, \rho_{\Theta}^{-1}\right)
\right)(x_{\Theta} \omega) \; ,
\end{equation}
is independent of the choice of $x_{\Theta } \in b(\Theta, \Gamma,
\omega)$.
\begin{fleuveIII}
\item In this part we suppose that
\begin{equation}\label{P5-5}
\lambda^i(\Omega_i) \;=\; 1 \quad \text{for every } i\in \mathbb{Z}^d\;.
\end{equation}
[As remarked below, this is not a big loss of generality.]
Let $\mathcal{N}$ be the set of probability measures $\mu$ on
$(\Omega,\mathcal{F})$ such that
\begin{equation}\label{N-set}
\mu\,\lambda_j\bigl[\comp{B(j,V)}\bigr]=0
\quad\text{for every } V\in\mathcal{S} \text{ and } j\in\comp{V}\;.
\end{equation}
Then, within this class, consistency is equivalent to consistency with
singletons:
\begin{equation}\label{cons-mea}
\mathcal{N}\,\cap\, \mathcal{G}\bigl(\{\rho_\Lambda\,
\lambda_\Lambda\}_{\Lambda\in\mathcal{S}}\bigr)
\;=\; \mathcal{N}\,\cap\, \Bigl\{\mu\in\mathcal{P}(\Omega,\tribu) :
\mu\,(\rho_i\, \lambda_i) = \mu, \;
i\in\mathbb{Z}^d\Bigr\}.
\end{equation}
\end{fleuveIII}
\end{theorem}

\begin{remarks}
\begin{itemize}
\item[]
\item In particular, if $x_{\Theta } \in b(\Theta, \Gamma, \omega)$,
  formulas \reff{r.1}--\reff{r.2} yield
  \begin{equation}\label{gr.10}
  \rho_{\Theta \cup \Gamma}(x_{\Theta} \omega) \;=\;
  \frac{\rho_{\Gamma}(x_{\Theta} \omega)}
  {\lambda_{\Gamma}\left(\rho_{\Gamma} \,
    \rho_{\Theta}^{-1}\right)(x_{\Theta} \omega)}\;,
  \end{equation}
  a formula already present in Theorem (1.33) of Georgii (1988).
\item A simple recursive argument shows that the translation invariance
  of the measures $\lambda^i$ and the functions $\rho_i$ imply that of the
  functions $\rho_\Lambda$.
\item When $E$ is finite and each $\lambda^i$ is the counting measure,
  results {\rm (I)} and {\rm (II)} were obtained by Dachian and Nahapetian
  (2004).  In the strictly positive case (everybody is good) we
  recover the results of the appendix of Fern\'andez and Maillard
  (2004). See Section \ref{sec5} for more details.
\item As remarked by Georgii [1988, Remark (1.28) (3)]\nocite{geo88},
  the normalization condition \reff{P5-5} is equivalent to the existence of
  functions $r_i(\omega_i)>0$ with $0<\lambda^i(r_i)<\infty$. Indeed,
the
  definition $\widetilde\rho_i = \rho_i/r_i$
  leads to the identity $\rho_i\, \lambda_i = \widetilde\rho_i\,
   \widetilde\lambda_i$ with $\widetilde\lambda^i(\Omega_i) = 1$.
  Such functions $r_i$ exist, for instance, if the measures $\lambda^i$
  are $\sigma$-finite.
\item If $E$ is compact, then usually both the measures
  $\lambda^i$ and the functions $\rho_i$ are bounded. In such a situation
  the continuity of $\rho_\Lambda$ and   $h$ implies the continuity of
  $\left(\rho_\Lambda\lambda_\Lambda\right)(h\mid \cdot\,)$ and the
  specification $\left\{ \rho_{\Lambda}
    \lambda_{\Lambda}\right\}_{\Lambda \in \mathcal{S}}$ is a
  \emph{Feller specification}.  If $E$ is finite such specifications
  are also \emph{quasilocal}.  See van Enter, Fern\'andez and Sokal
  (1993) for a survey of these notions and their relation to
  Gibbsianness.
\item Preston (2004)\nocite{pres04}, in an unpublished preprint,
  proves a rather strong result related to our part {\rm (III)}.
  The author takes a reconstructive point of view
  ---~the functions $\rho_i$ come from a pre-existing specification~---
  and determines conditions under which consistency coincides with
  singleton consistency. His framework is more general than ours in that
  a product structure is not demanded. On the other hand, the hypothesis
  imposed by Preston to the consistent measure involves all kernels,
  and not only the singletons as in \reff{N-set}.
  \end{itemize}
\end{remarks}

The following proposition explains in which sense the
order-consistency condition {\rm (H2)} is natural for a
specification satisfying {\rm (H1)}. Indeed, if a specification
$\gamma=(\rho_\Lambda\, \lambda_\Lambda)_{\Lambda\in\mathcal{S}}$
is such that the family $\{\rho_i\}_{i\in\mathbb{Z}^d}$ satisfies
{\rm (H1)} and the good configurations have a positive probability
and, then an almost sure version of {\rm (H2)} is fulfilled. In
particular, when $E$ is countable, {\rm (H2)} is fully satisfied.
\begin{proposition}\label{prop1}
Assume that $\lambda^i(b(i,\{j\},\alpha))>0$ for all
$\alpha\in\Omega$ and $i,j\in\mathbb{Z}^d$ such that $i\neq j$.
Then, for $\lambda_{\{i,j\}}(\,\cdot\mid\alpha)$-almost all
$\omega\in\Omega$
\begin{equation}
\frac{\rho_i(\omega)\,
\rho_j(x_i\omega)}{\rho_i(x_i\omega)\,\lambda_j(\rho_j\,\rho_i^{-1})(x_i\omega)}
\;=\;\frac{\rho_j(\omega)\,
\rho_i(x_j\omega)}{\rho_j(x_j\omega)\,\lambda_i(\rho_i\,\rho_j^{-1})(x_j\omega)}
\end{equation}
for all $x_i\in b(i,\{j\},\omega)$ and $x_j\in
b(j,\{i\},\omega)$. In particular, when $E$ is countable and each
$\lambda^i$ is the counting measure, {\rm (H2)} is satisfied for
all $\omega\in\Omega$.
\end{proposition}

As an illustration of our results, we present a family of singletons
satisfying the hypotheses of Theorem \ref{th1} but not fitting any of
the existing (re)construction schemes.  The main value of this example is to
provide a concrete manifestation of the different hypotheses of the theorem.

\begin{example}\label{ex1}
  Let $E=[0,1]$ and $\mathcal{E}$ be its Borel $\sigma$-algebra.  For
  each $i\in\zzz$ we take $\lambda^i$ equal to the Lebesgue measure
  and define
\begin{equation}
\rho_{i}(\omega) \;=\; \begin{cases}
2 \ind{[0,1/2]}(\omega_i) & \text{if } \left|\left\{j : \omega_j >
  1/2\right\}\right| = \infty \;,\\
2 \ind{]1/2,1]}(\omega_i) & \text{otherwise}\;.
\end{cases}
\end{equation}
Let us see that these functions satisfy the hypotheses of Theorem
\ref{th1}. The measurability of each $\rho_i$ and the normalization
\reff{spe33} are readily verified.  We check {\rm (H1)} and {\rm (H2)} for
$\omega$ such that $\left|\left\{j : \omega_j > 1/2\right\}\right| =
\infty$, the complementary case is analogous.  For such $\omega$ we
see that for all $j \in \zzz^{d}$ and $V \in
\mathcal{S}\left(\comp{\{j\}}\right)$:
\begin{itemize}
\item[(i)] $\rho_j(\sigma_V \omega)>0$ for all $\sigma_V \in
\Omega_V$ if and only if $\omega_j \in [0,1/2]$.
\item[(ii)] If $\omega_j \in [0,1/2]$ then
$\lambda_i\left(\rho_i \rho_{j}^{-1}\right)(\sigma_V \omega) = 1/2$
for all $\sigma_V \in \Omega_V$.
\end{itemize}
It follows that {\rm (H1)} is verified with $b(j,V,\omega) = [0,1/2]$.
Furthermore,
\begin{equation}\label{ex4}
R_{i}^{j}(\omega) \;=\;  \left\{\begin{array}{ll}
1/2 & \mbox{if } \omega_j \in [0,1/2]\\
\infty & \mbox{otherwise}\;,
\end{array}\right.
\end{equation}
satisfying \rm{ (H2)(b)}, and, if $x_i, x_j \in
[0,1/2]$,
\begin{equation}\label{ex3}
\begin{aligned}
\frac{\rho_{i}(\omega) \, \rho_j (x_i \omega)}{\rho_i (x_i
\omega)\, \lambda_{j}\left(\rho_{j} \, \rho_{i}^{-1}\right)(x_i
\omega)} &\;=\; \frac{\rho_{j}(\omega) \, \rho_i (x_j \omega)}{\rho_j
(x_j \omega)\, \lambda_{i}\left(\rho_{i} \, \rho_{j}^{-1}\right)(x_j
\omega)}\nonumber\\[5pt]
&\;=\; \left\{\begin{array}{ll}
4 & \mbox{if } \omega_i,\omega_j \in [0,1/2]\\
0 & \mbox{otherwise}\;,
\end{array}\right.
\end{aligned}
\end{equation}
in agreement with hypotheses \rm{(H2)(a)}.
We observe that our construction indeed leads to
\begin{equation}
\rho_{\Lambda}(\omega) \;=\; \begin{cases}
2^{|\Lambda|} \ind{[0,1/2]^\Lambda}(\omega_\Lambda) & \text{if }
\left|\left\{j : \omega_j >
  1/2\right\}\right| = \infty \;,\\
2^{|\Lambda|} \ind{]1/2,1]^\Lambda}(\omega_\Lambda) & \text{otherwise}\;.
\end{cases}
\end{equation}
\end{example}

\section{Comparison with previous results}\label{sec5}
Our results are the generalization of those obtained by Dachian and Nahapetian
(2004)\nocite{dacnah04}
in the case of finite single-site space $E$ and each $\lambda^i$ equals to the
counting measure. In this
framework our hypotheses {\rm (H1)} reduces to the positivity requirement
\reff{spe16}. A family of single-site
weights $\{\rho_i\}_{i\in\mathbb{Z}^d}$ satisfying such a property is termed
\emph{very weakly positive} by
the authors. In our notation, their result (obtained by combining their
Proposition 18, Theorem 19 and Theorem 21)
is the following
\begin{proposition}[Dachian and Nahapetian (2004)]\label{DH}
Let $\{\rho_i\}_{i\in\mathbb{Z}^d}$ be a family of very weakly positive
probability weights which are normalized
in the sense that
\begin{equation}\label{norm}
\lambda^i(\rho_i)\equiv 1\;.
\end{equation}
Then:
\begin{fleuveII}
\item There exists a unique family $\left\{ \rho_{\Lambda}
\right\}_{\Lambda   \in \mathcal{S}}$ of measurable functions
$\rho_\Lambda:\Omega\to
[0,\infty[$, with $\rho_{\{i\}}=\rho_i$, such that the family of kernels
$\left\{ \rho_{\Lambda}
\lambda_{\Lambda}\right\}_{\Lambda \in \mathcal{S}}$ is a specification if and
only if
\begin{equation}\label{DH1}
\begin{aligned}
&\rho_i(x_j\, u_i\, \omega)\, \rho_j(u_j\, u_i\, \omega)\, \rho_i(x_i\,
u_j\, \omega)\, \rho_j(x_j\, x_i\, \omega)
\nonumber\\
&\qquad\qquad\;=\;
\rho_j(x_i\, u_j\, \omega)\, \rho_i(u_i\, u_j\, \omega)\, \rho_j(x_j\, u_i\,
\omega)\, \rho_i(x_i\, x_j\, \omega)
\;,
\end{aligned}
\end{equation}
for every $i,j\in\mathbb{Z}^d$, $u_i\in\Omega_i$, $u_j\in\Omega_j$, $x_i\in
b(i,\{j\},\omega)$ and
$x_j\in b(j,\{i\},\omega)$.
\item The functions $\left\{ \rho_{\Lambda}
\right\}_{\Lambda   \in \mathcal{S}}$ are continuous if and only if the
functions $\left\{ \rho_{i}
\right\}_{i   \in \mathbb{Z}^d}$ are.
\end{fleuveII}
\end{proposition}
For finite state-space $E$, we claim that this proposition coincides with parts
{\rm (I)} and {\rm (II)}
of Theorem \ref{th1}.
To prove this, it suffices to show that \reff{DH1} is equivalent to our
hypotheses {\rm (H2)}. In fact, the
``only if" part of Proposition \ref{DH} implies that \reff{DH1} is
satisfied whenever {\rm (H2)} is. Thus, we only
need to show that \reff{DH1} implies our condition {\rm (H2)}. This is easily
seen. Indeed, if we divide both sides
of \reff{DH1} by $\rho_i(x_i\, u_j\, \omega)\, \rho_j(x_j\, u_i\, \omega)$ and
sum them over $u_i\in\Omega_i$ and
$u_j\in\Omega_j$, we obtain thanks to the normalization \reff{norm},
\begin{equation}
\rho_i(x_i\, x_j\, \omega)\, \sum_{u_j\in\Omega_j}\frac{\rho_j(x_i\, u_j\,
\omega)}{\rho_i(x_i\, u_j\, \omega)}
\;=\;
\rho_j(x_j\, x_i\, \omega)\, \sum_{u_i\in\Omega_i}\frac{\rho_i(x_j\, u_i\,
\omega)}{\rho_i(x_i\, u_j\, \omega)}\;,
\end{equation}
that is
\begin{equation}\label{DH2}
\rho_i(x_i\, x_j\, \omega)\, \lambda_j(\rho_j\, \rho_i^{-1})(x_i\, \omega)
\;=\;
\rho_j(x_j\, x_i\, \omega)\, \lambda_i(\rho_i\, \rho_j^{-1})(x_j\, \omega)\;.
\end{equation}
Dividing term-by-term \reff{DH1} by \reff{DH2}, we arrive to \reff{r.20}.

\bigskip
A related but complementary result is contained in Appendix A of
Fer\-n\'an\-dez Maillard (2004). In this appendix, the
configuration space $\Omega$ is an arbitrary subset of
$E^{\mathbb{Z}^d}$, for a general measurable space $E$. Thus,
$\Omega$ may describe local exclusion rules or grammars. Kernels
are supposed to be strictly positive on the whole of $\Omega$ and
therefore our hypotheses {\rm (H2)} becomes
\begin{equation}\label{gibbs31}
\frac{\rho_{i}}{\lambda_{i} \left( \rho_{i} \,
  \rho_{j}^{-1}\right)}(\omega)
\;= \;\frac{\rho_{j}}{ \lambda_{j} \left( \rho_{j} \,
  \rho_{i}^{-1}\right)}(\omega)\;,
\end{equation}
for every $i,j$ in   $\mathbb{Z}^{d}$ and every $\omega\in\Omega$. The result
is

\begin{proposition}[Fern\'andez and Maillard (2004)]\label{thap}
Let $\{\rho_i\}_{i\in\mathbb{Z}^d}$ be a family of measurable functions which
are normalized ---
$\lambda^i(\rho_i)\equiv 1$ --- and satisfy \reff{gibbs31} and the following
bounded-positivity properties.
For every $i,j \in \mathbb{Z}^{d}$,
\begin{equation}\label{gibbs27}
 \inf_{\omega \in \Omega}
\lambda_{j} \left( \rho_{j} \, \rho_{i}^{-1} \right)(\omega)\;>\;0\;,
\end{equation}
and
\begin{equation}\label{gibbs27.1}
\sup_{\omega \in \Omega} \lambda_{j} \left( \rho_{j} \,
\rho_{i}^{-1} \right)(\omega) \;<\; +\infty\;.
\end{equation}
Then:
\begin{fleuveII}
\item There exists a unique family $\left\{ \rho_{\Lambda}
\right\}_{\Lambda   \in \mathcal{S}}$ of measurable functions
$\rho_\Lambda:\Omega\to
[0,\infty[$, with $\rho_{\{i\}}=\rho_i$, such that the family of kernels
$\left\{ \rho_{\Lambda}
\lambda_{\Lambda}\right\}_{\Lambda \in \mathcal{S}}$ is a specification.
\item  If the functions $\rho_i$ are continuous and $\int
    \sup_\omega (\rho_i\,\rho_j^{-1})(\sigma_i\omega_{\comp{\{i\}}})
    \,\lambda^i(d\sigma_i) \;<\:\infty$ for all $i,j\in\mathbb{Z}^d$,
    then the functions $\rho_\Lambda$, and thus the specification
    $\gamma$, are continuous.
\item $\mathcal{G}(\left\{ \rho_{\Lambda}
\lambda_{\Lambda}\right\}_{\Lambda \in \mathcal{S}}) = \left\{ \mu \in
\mathcal{P}(\Omega,
      \tribu) : \mu (\rho_i\, \lambda_i) = \mu \mbox{ for all } i \in
      \mathbb{Z}^{d} \right\}$.
\item For each $\Lambda\in\mathcal{S} $ there exist constants
    $C_\Lambda, D_\Lambda >0$ such that $C_\Lambda\,\rho_k(\omega) \le
    \rho_\Lambda(\omega) \le D_\Lambda\,\rho_k(\omega)$ for all
    $k\in\Lambda$ and all $\omega\in\Omega$.
\end{fleuveII}
\end{proposition}
If $\Omega$ is the full product space $E^{\mathbb{Z}^d}$ and every configuration
is allowed
---~$b(j,V,\omega) = \Omega_j$ $\forall\; V\in \mathcal{S}$, $j\in \comp{V}$,
$\omega\in\Omega$~---
Proposition \ref{thap} and Theorem \ref{th1} coincide. But, otherwise, these two
results have different
ranges of application. Indeed, models with local exclusion rules (for instance,
no two nearest neighbor
simultaneously occupied) are covered by Proposition \ref{thap}, but do not
satisfy the hypotheses of Theorem
\ref{th1}. The reason is that each $b(j,V,\omega)$ is empty if $V\neq \emptyset$
or if $\omega$ violates the
exclusion rules. On the other hand, models with ``asymptotic" exclusion
rules, like in Example \ref{ex1},
fall outside the scope of Proposition \ref{thap}.

\section{Proof of Theorem \ref{th1}}

We need tree lemmas to build the proof of our theorem. We start with the
crucial one showing that the algorithm \reff{r.1}--\reff{r.2}
recursively leads to multi-site generalizations of hypotheses (H1) and
(H2).

\begin{lemma}\label{spelem1}
Let $\{\rho_i\}_{i \in \zzz^{d}}$ be a family of $\tribu$-measurable
  functions $\rho_i:\Omega \to [0,\infty[$ satisfying hypotheses {\rm (H1)}
  and {\rm (H2)}.  Then

\begin{fleuve11}
\item The equations
\begin{equation}\label{spe27}
\rho_{\Lambda \cup \{i\}}(\omega) \;=\;
\frac{\rho_{\Lambda}(\omega)}{R_{\Lambda}^{i}(\omega)}\;,
\end{equation}
\begin{equation}\label{spe28}
R_{\Lambda}^{i}(\omega) \;=\;
\left(\frac{\rho_{\Lambda}}{\rho_{i}} \times \lambda_{i}
  \left(\rho_{i} \, \rho_{\Lambda}^{-1}\right)\right)
\left(x_{\Lambda} \omega \right)\;,
\end{equation}
$i\not\in\Lambda$, recursively define for each $\Lambda\in\mathcal{S}$
measurable functions $\rho_{\Lambda} : \Omega \rightarrow [0,+\infty[$
and $R_{\Lambda}^{i} : \Omega \longrightarrow ]0,+\infty]$, the latter
being independent of the choice of $x_{\Lambda} \in b(\Lambda, \{i\},
\omega)$.
\item The functions defined above satisfy that for each $\omega \in
  \Omega$, $V \in \mathcal{S}(\comp{\Lambda})$, $x_{\Lambda} \in
  b(\Lambda, V, \omega)$,
\begin{equation}\label{spe22}
\rho_{\Lambda}\left(x_{\Lambda} \sigma_V \omega\right)>0 \quad
\forall \, \sigma_V \in \Omega_V\;,
\end{equation}
 and, for each $j \in \Lambda$ and $i\in\zzz^d$, $i\neq j$,
\begin{equation}\label{spe23}
\inf \left\{\lambda_{i}\left(\rho_{i} \, \rho_{j}^{-1} \,
    \lambda_{j}\left(\rho_{j} \,
\rho_{\Lambda_{j}^{*}}^{-1}\right)\right)\left(x_{\Lambda} \sigma_{V}
\omega \right) : \sigma_{V} \in \Omega_{V} \right\}>0
\end{equation}
and
\begin{equation}\label{spe24}
\sup \left\{\lambda_{i}\left(\rho_{i} \, \rho_{j}^{-1} \,
    \lambda_{j}\left(\rho_{j}
\, \rho_{\Lambda_{j}^{*}}^{-1}\right)\right)\left(x_{\Lambda}
\sigma_{V} \omega \right) : \sigma_{V} \in
\Omega_{V} \right\}<\infty\;,
\end{equation}
with the convention that $\rho_{\emptyset} \equiv 1$.
\item More generally,
\begin{equation}\label{r.51}
\rho_{\Theta \cup \Gamma}(\omega) \;=\;
\frac{\rho_{\Theta}(\omega)}{R_{\Theta}^{\Gamma}(\omega)}
\end{equation}
with
\begin{equation}\label{r.50}
R_{\Theta}^{\Gamma}(\omega) \;=\;
\left(\frac{\rho_\Theta}{\rho_\Gamma} \times
\lambda_{\Gamma}\left(\rho_{\Gamma} \, \rho_{\Theta}^{-1}\right)
\right)(x_{\Theta} \omega)\;,
\end{equation}
for each $\Theta \in \mathcal{S}$, $\Gamma \in
\mathcal{S}(\comp{\Theta})$ and $\omega \in \Omega$.  The RHS of
\reff{r.50} is independent of $x_\Theta\in b(\Theta,\Gamma,\omega)$.
\end{fleuve11}
\end{lemma}

\prf We will prove the Lemma by induction over $|\Lambda| \geq 1$.
In (3) we assume $\Theta\cup\Gamma= \Lambda\cup \{i\}$ for some
$i\not\in\Lambda$.  Note that, in particular, (3) implies that the
value of the functions $\rho_\Lambda$ do not depend on the order in
which the sites of $\Lambda$ are swept during the recursive
construction.

The initial inductive step is immediate: If $\Lambda = \{j\}$, item
(1) amounts to the identity \reff{r.31} (with $i\leftrightarrow j$)
and \reff{spe28} is just the definition of $R_j^i$.  Item (2)
coincides with hypothesis (H1) while item (3) is the identity
\reff{c.2} which remains valid even if some numerator is zero or some
denominator is infinity.

Suppose now (1)--(3) valid for all finite subsets of $\zzz^d$
involving up to $n$ sites.  Consider $\Lambda \in \mathcal{S}$ of
cardinality $n+1$, $i\not\in\Lambda$, $V \in \mathcal{S}$ with $V
\subset \comp{\Lambda}$ and some $x_{\Lambda} \in b(\Lambda, V ,
\omega)$.  We observe that, by the very definition of $b(\Lambda, V ,
\omega)$ [see \reff{spe13}],
\begin{equation}\label{co.1}
x_j\in b(j, \Lambda^*_j \cup V, \omega) \quad\mbox{and}\quad
x_{\Lambda^*_j}\in b(\Lambda^*_j , V\cup\{j\}, \omega)
\end{equation}
for each site $j \in \Lambda$.  The leftmost statement implies, by
hypothesis (H1), that
\begin{equation}\label{spe46}
\rho_{j}(x_{j} \sigma_{\Lambda^*_j \cup V} \omega
)>0, \quad \forall \, \sigma_{\Lambda^*_j \cup V} \in
\Omega_{\Lambda^*_j \cup V}\;,
\end{equation}
and, if $i\neq j$,
\begin{equation}\label{spe46bis}
\left.\begin{array}{l} \inf  \\ \sup \end{array} \right\}
\Bigl\{\lambda_{i}\left(\rho_{i} \, \rho_{j}^{-1}\right)
(x_{j} \sigma_{\Lambda^*_j \cup V} \omega ) :
\sigma_{\Lambda^*_j \cup V} \in \Omega_{\Lambda^*_j \cup V} \Bigr\}
\left\{\begin{array}{l} >0  \\ <\infty\;. \end{array} \right.
\end{equation}
On the other hand, the rightmost statement in \reff{co.1} and the
inductive hypothesis (2) imply that
\begin{equation}\label{spe48}
\left.\begin{array}{l} \inf  \\ \sup \end{array} \right\}
\Bigl\{ \lambda_{j}\left(\rho_{j} \, \rho_{\Lambda^*_j}^{-1}\right)
(x_{\Lambda^*_j} \sigma_{V \cup \{i\}} \omega ) :
\sigma_{V \cup \{i\}} \in \Omega_{V \cup \{i\}} \Bigr\}
\left\{\begin{array}{l} >0  \\ <\infty\;. \end{array} \right.
\end{equation}

\paragraph{Proof of (2)~:\,\,}
Combining \reff{spe46} and \reff{spe48} we see that the quotient
\begin{equation}\label{spe49}
\rho_{\Lambda}\left(x_{\Lambda} \sigma_{V} \omega\right)\;\triangleq\;
\frac{\rho_{j}\left(x_{\Lambda} \sigma_{V} \omega \right)}
{\lambda_{j}\left(\rho_{j} \,\rho_{\Lambda^*_j}^{-1}\right)
\left(x_{\Lambda} \sigma_{V} \omega \right)}
\end{equation}
satisfies
\begin{equation}\label{a.a}
0\;<\; \rho_{\Lambda}\left(x_{\Lambda} \sigma_{V} \omega\right)
\;<\; \infty\;,
\end{equation}
while \reff{spe46} and \reff{spe48} imply that
\begin{equation}\label{spe49bis}
\left.\begin{array}{l} \inf  \\ \sup \end{array} \right\}
\Bigl\{ \lambda_{i}\left(\rho_{i} \,\rho_{j}^{-1} \,
\lambda_{j}(\rho_{j} \, \rho_{\Lambda^*_j}^{-1})\right)
\left(x_{\Lambda} \sigma_{V} \omega \right) :
\sigma_{V} \in \Omega_{V} \Bigr\}
\left\{\begin{array}{l} >0  \\ <\infty\;. \end{array} \right.
\end{equation}
Together \reff{spe49} and \reff{spe49bis} yield
\begin{equation}\label{ro.10}
\left.\begin{array}{l} \inf  \\ \sup \end{array} \right\}
\Bigl\{ \lambda_i\left(\rho_i \, \rho_\Lambda^{-1}\right)
\left(x_{\Lambda} \sigma_{V} \omega \right) :
\sigma_{V} \in \Omega_{V} \Bigr\}
\left\{\begin{array}{l} >0  \\ <\infty\;. \end{array} \right.
\end{equation}

\paragraph{Proof of (1)~:\,\,}
We consider now $V=\{i\}$ with $i\not\in \Lambda$.  Inequalities
\reff{ro.10} and the symmetry relation \reff{a.a} imply that
$\left(\rho_{\Lambda}\lambda_i(\rho \, \rho_\Lambda^{-1})\right)
\left(x_{\Lambda} \sigma_i \omega \right)>0$ for all $\sigma_i\in
b(\Lambda, \{i\} ,\omega)$ and thus it makes sense to define
\begin{equation}\label{spe51}
R_{\Lambda}^{i}(\omega) \;=\; \left(\frac{\rho_{\Lambda}}{\rho_{i}}
  \times \lambda_{i} \left(\rho_{i} \,
    \rho_{\Lambda}^{-1}\right)\right)
\left(x_{\Lambda} \omega\right)
\end{equation}
which may be infinite but, due to \reff{a.a} and \reff{ro.10}, is
never zero.  We conclude that the function $\rho_{\Lambda\cup\{i\}}$
defined by \reff{spe27} takes values on $[0,\infty[$.

We must prove that definition \reff{spe51} is indeed independent of
the choice of $x_\Lambda\in b(\Lambda, \{i\} , \omega)$.  We analyze
first the case $R_{\Lambda}^{i}(\omega) <\infty$.  For each $j \in
\Lambda$ and each $\sigma_{i} \in \Omega_{i}$ we have, by the
inductive hypothesis (1),
\begin{equation}\label{spe58}
\rho_{\Lambda}(x_{\Lambda} \sigma_{i} \omega ) \;=\;
\frac{\rho_{j}}{\lambda_{j}\left(\rho_{j} \,
    \rho_{\Lambda_{j}^{*}}^{-1}\right)}
(x_{\Lambda} \sigma_{i} \omega)\;.
\end{equation}
Furthermore, combining \reff{spe51} and \reff{spe58} we obtain
\begin{equation}
R_{\Lambda}^{i}(\omega) \;=\; \left(\frac{\rho_{j}}{\rho_{i} \times
    \lambda_{j}\left(\rho_{j} \,
\rho_{\Lambda_{j}^{*}}^{-1}\right)} \times \lambda_{i}
\left(\frac{\rho_{i} \times \lambda_{j}\left(\rho_{j} \,
\rho_{\Lambda_{j}^{*}}^{-1}\right)}{\rho_{j}}\right)\right)
(x_{\Lambda}\omega )\;.
\end{equation}
We now use \reff{c.2}, namely
$\rho_i/R_i^j = \rho_j/R_j^i$, and make use of the
$\tribu_{\comp{\{i\}}}$-measurability of $R_{i}^{j}$ to pass it through
the $\lambda_i$-integration.  We get
\begin{equation}\label{gr.300}
\begin{aligned}
R_{\Lambda}^{i}(\omega)
&\;=\;
\left(\frac{R_{j}^{i}}{\lambda_{j}\left(\rho_{j} \,
      \rho_{\Lambda_{j}^{*}}^{-1}\right)} \times \lambda_{i}
\left(\frac{\lambda_{j}\left(\rho_{j} \,
      \rho_{\Lambda_{j}^{*}}^{-1}\right)}{R_{j}^{i}}\right)\right)
(x_{\Lambda}\omega )\nonumber\\
&\;=\; \left(\frac{R_{j}^{i}}{\lambda_{j}\left(\rho_{j} \,
      \rho_{\Lambda_{j}^{*}}^{-1}\right)} \times \lambda_{\{i,j\}}
\left(\frac{\rho_{j} \,    \rho_{\Lambda_{j}^{*}}^{-1}}
{R_{j}^{i}}\right)\right) (x_{\Lambda}\omega )\;.
\end{aligned}
\end{equation}
In the last equality we used the factorization property
\reff{spe15} of the free kernel and the
$\tribu_{\comp{j}}$-measurability of $R_j^i$.  The final expression is
manifestly independent of the actual value of $x_j$.  Since $j$ is an
arbitrary site of $\Lambda$, we conclude that $R_{\Lambda}^{i}$ is
$\tribu_{\comp{\Lambda}}$-measurable.

Let us turn now to the case $R_\Lambda^i(\omega)=\infty$.  This
happens if, and only if, $\rho_i(x_\Lambda\omega)=0$.  We must prove
that, in this case, $\rho_i(\widetilde x_j x_{\Lambda^*_j}\omega)=0$
for any $j\in\Lambda$ and any $\widetilde x_j\in\Omega_j$ such that
$\widetilde x_j x_{\Lambda^*_j}\in b(\Lambda, \{i\} , \omega)$.  But,
by the definition of $b(\Lambda, \{i\} , \omega)$, for every
$j\in\Lambda$
\begin{equation}\label{ro.a}
R_i^j(x_\Lambda\omega) < \infty \quad\mbox{and}\quad
\rho_j(x_\Lambda\omega) >0
\end{equation}
and
\begin{equation}\label{ro.b}
R_i^j(\widetilde x_j x_{\Lambda^*_j}\omega) < \infty
\quad\mbox{and}\quad
\rho_j(\widetilde x_j x_{\Lambda^*_j}\omega) >0\;.
\end{equation}
We can now establish the following chain of implications:
\begin{equation}\label{gr.301}
\begin{aligned}
\lefteqn{\rho_i(x_\Lambda\omega)=0 \quad\Longrightarrow\quad
R_j^i(x_\Lambda \omega)=\infty}\nonumber\\
&&  \quad\Longrightarrow\quad
R_j^i(\widetilde x_j x_{\Lambda^*_j}\omega) =\infty
\quad\Longrightarrow\quad \rho_i(\widetilde x_j x_{\Lambda^*_j}\omega) =0\;.
\end{aligned}
\end{equation}
The first implication results from \reff{ro.a} and the symmetry
relation \reff{c.2}, the second one is a consequence of the
$\tribu_{\comp{\{j\}}}$-measurability of $R_j^i$ and the last one
follows from \reff{ro.b} and \reff{c.2}.

\paragraph{Proof of (3)~:\,\,}
We consider $\Theta$ and $\Gamma$ disjoint, non-empty, with
$\left|\Theta\cup\Gamma\right|=n+1$, for $n\ge 2$ (the case $n=1$ was
analyzed at the begining).  We have to prove that if
$\Theta\cup\Gamma=\widetilde\Theta\cup\widetilde\Gamma$ with
$\widetilde\Theta$ and $\widetilde\Gamma$ disjoint, then
\begin{equation}
\frac{\rho_\Theta}{R_\Theta^\Gamma} \;=\;
\frac{\rho_{\widetilde\Theta}}{R_{\widetilde\Theta}^{\widetilde\Gamma}}
\;.
\end{equation}
As the argument is symmetric in $\Theta$ and $\Gamma$ we can assume
that $\left|\Theta\right|\ge2$, in which case, modulo iteration, it
is enough to prove that for $k\in\Theta$
\begin{equation}\label{gr.2}
\frac{\rho_\Theta}{R_\Theta^\Gamma} \;=\;
 \frac{\rho_{\Theta_k^*}}{R_{\Theta_k^*}^{\Gamma\cup\{k\}}} \;.
\end{equation}
Let us fix some $\omega\in\Omega$ and $x_\Theta\in
b(\Theta,\Gamma,\omega)$.  The inductive definition
\reff{r.51}--\reff{r.50} immediately yields the identity
\begin{equation}\label{gr.3}
\frac{\rho_\Theta(\omega)}{\rho_\Theta(x_\Theta\omega)} \;=\;
\frac{\rho_{\Theta^*_k}(\omega)\rho_k(x_{\Theta^*_k}\omega)}
{\rho_{\Theta^*_k}(x_{\Theta^*_k}\omega)\rho_k(x_\Theta\omega)} \;.
\end{equation}
In addition we need the following identity
\begin{equation}\label{gr.5}
\lambda_\Gamma\left(\rho_\Gamma\rho_\Theta^{-1}\right)
(x_\Theta\omega)\;=\;
\lambda_\Gamma\left(\rho_\Gamma\rho_k^{-1}\right)
(x_\Theta\omega)\,
\lambda_{\Gamma\cup\{k\}}\left(\rho_{\Gamma\cup\{k\}}
\rho_{\Theta_k^*}^{-1}\right) (x_{\Theta_k^*}\omega)\;.
\end{equation}
This is proved as follows.  We start from the relation
\begin{equation}\label{gr.19.0}
\left(\rho_\Gamma \, \rho_{\Theta}^{-1} \right)
(x_{\Theta} \omega ) \;=\;
\left(\rho_\Gamma \, \rho_{k}^{-1} \,
\lambda_{k}\left( \rho_{k} \, \rho_{\Theta_{k}^{*}}^{-1} \right)
\right)(x_{\Theta} \, \omega )
\end{equation}
which is an immediate consequence of the inductive hypotheses
\reff{r.51}--\reff{r.50} [see \reff{gr.10}].  As $x_k\in
b(k,\Gamma\cup\Theta^*_k,\omega)$, the LHS
is well defined for every $\omega_\Gamma$.  We can, therefore,
integrate both sides and conclude that
\begin{equation}\label{gr.19}
\lambda_\Gamma\left(\rho_\Gamma \, \rho_{\Theta}^{-1} \right)
(x_{\Theta} \omega ) \;=\;
\lambda_\Gamma\left(\rho_\Gamma \, \rho_{k}^{-1} \,
\lambda_{k}\left( \rho_{k} \, \rho_{\Theta_{k}^{*}}^{-1} \right)
\right)(x_{\Theta} \, \omega )\;.
\end{equation}
Next we observe that
\begin{equation}\label{gr.20}
\frac{\rho_\Gamma(x_\Theta\sigma_\Gamma\omega)}
{\rho_k(x_\Theta\sigma_\Gamma\omega)} \;=\;
\frac{\lambda_\Gamma\left(\rho_\Gamma\,\rho_k^{-1}\right)
(x_\Theta\omega)}{R_k^\Gamma(x_\Theta\sigma_\Gamma\omega)}
\end{equation}
for all $\sigma_\Gamma\in\Omega_\Gamma$.  Again, this is a consequence
of the inductive validity of \reff{r.51}--\reff{r.50} which, in particular,
also implies that if $x_\Theta\in b(\Theta,\Gamma,\omega)$,
\begin{equation}
R_\Gamma^k (x_\Theta\omega) \;=\;
\lambda_\Gamma\left(\rho_\Gamma\,\rho_k^{-1}\right)(x_\Theta\omega)\;.
\end{equation}
To obtain \reff{gr.5} we must insert \reff{gr.20} into \reff{gr.19}
and use that by \reff{spe15}
$\lambda_\Gamma\lambda_k=\lambda_{\Gamma\cup\{k\}}$.
\smallskip

The combination of \reff{gr.3} and \reff{gr.5} yields, thanks to the
inductive definition of $\rho_{\Gamma\cup\{k\}}(x_{\Theta_k^*}\omega)$,
\begin{equation}\label{gr.30}
\frac{\rho_\Theta}{R_\Theta^\Gamma} (\omega)\;=\;
\frac{\rho_{\Theta^*_k}(\omega)\,\rho_{\Gamma\cup\{k\}}(x_{\Theta^*_k}\omega)}
{\rho_{\Theta^*_k}(x_{\Theta^*_k}\omega)\,
\lambda_{\Gamma\cup\{k\}}\left(\rho_{\Gamma\cup\{k\}}
\rho_{\Theta_k^*}^{-1}\right) (x_{\Theta_k^*}\omega)}\;.
\end{equation}
Due to the inductive definition \reff{r.50} of
$R_{\Theta^*_k}^{\Gamma\cup\{k\}}$, the RHS of \reff{gr.30} is
precisely the RHS of \reff{gr.2}.  This concludes the proof of (3), at
least when $R_\Theta^\Gamma(\omega)<\infty$.  But in fact the argument
leading to identity \reff{gr.30} remains valid also when
$R_\Theta^\Gamma(\omega)$ is infinite. In this case we have the
following chain of implications:
\begin{equation}\label{gr.40}
R_\Theta^\Gamma(\omega)=\infty \quad\Longrightarrow\quad
\rho_{\Gamma\cup\{k\}}(x_{\Theta^*_k}\omega)\;=\;0
 \quad\Longrightarrow\quad
R_{\Theta^*_k}^{\Gamma\cup\{k\}}(\omega)=\infty\;.
\end{equation}
The first implication is due to \reff{gr.30} while the second one
follows from the inductive definition of
$R_{\Theta^*_k}^{\Gamma\cup\{k\}}$. Display \reff{gr.40} proves
\reff{gr.2} when $R_\Theta^\Gamma(\omega)=\infty$.

The proof that $R_\Theta^\Gamma(x_\Theta\omega)$ is independent of
$x_\Theta\in b(\Theta,\Gamma,\omega)$ is completely analogous to the
preceding proof of (1).  We leave to the reader the pleasure of
obtaining a formula similar to \reff{gr.300} and a chain of
implications similar to \reff{gr.301} but changing $\Lambda\to\Theta$
and $i\to\Gamma$.  $\quad \Box$

\medskip
The following is a rather elementary property of conditional expectations.

\begin{lemma}\label{lem1}
  Let $\left\{\gamma_{\Lambda} \right\}_{\Lambda \in \mathcal{S}}$ be
  a specification, then for each $\Lambda \in \mathcal{S}, \, \Gamma
  \in \mathcal{S}_{\comp{\Lambda}}$ and bounded measurable functions $f,
  \, g$,
  \begin{equation} \label{gr.1}
\gamma_{\Lambda \cup \Gamma}\Bigr[f \, \gamma_{\Lambda}\bigr(
    \gamma_{\Gamma}(g)\bigr)\Bigr] \;=\; \gamma_{\Lambda \cup
      \Gamma}\Bigr[ g \,
    \gamma_{\Gamma}\bigr(\gamma_{\Lambda}(f)\bigr)\Bigr]\;.
    \end{equation}
\end{lemma}

\prf By the consistency of the specification and the
$\tribu_{\comp{\Lambda}}$-measurability of $\gamma_{\Lambda}\bigr(
\gamma_{\Gamma}(g)\bigr)$, we have
\begin{equation}\label{ro.5}
\begin{aligned}
\gamma_{\Lambda \cup \Gamma}\Bigr[f \,
\gamma_{\Lambda}\bigr(\gamma_{\Gamma}(g)\bigr)\Bigr]
&\;=\; \gamma_{\Lambda \cup \Gamma}\Bigr[\gamma_{\Lambda} \left( f \,
\gamma_{\Lambda}\bigr( \gamma_{\Gamma}(g)\bigr)\right)\Bigr]\nonumber\\
&\;=\; \gamma_{\Lambda \cup \Gamma}\Bigr[\gamma_{\Lambda}(f) \,
\gamma_{\Lambda}\Bigr(\gamma_{\Gamma}(g)\Bigr)\Bigr]\nonumber
\end{aligned}
\end{equation}
Similarly, the $\tribu_{\comp{\Lambda}}$-measurability of
$\gamma_{\Lambda}(f)$ and the consistency of the specification give
\begin{equation}
\begin{aligned}
\gamma_{\Lambda \cup \Gamma}\Bigr[\gamma_\Lambda(f) \,
\gamma_{\Lambda}\Bigr(\gamma_{\Gamma}(g)\Bigr)\Bigr]
&\;=\; \gamma_{\Lambda \cup
  \Gamma}\Bigr[\gamma_{\Lambda}\Bigr(\gamma_{\Lambda}(f) \,
\gamma_{\Gamma}(g)\Bigr)\Bigr]
\nonumber\\
&\;=\;\gamma_{\Lambda \cup \Gamma}\Bigr[\gamma_{\Lambda}(f) \,
\gamma_{\Gamma}(g)\Bigr]\;. \nonumber
\end{aligned}
\end{equation}
Identity \reff{gr.1} follows from the $f \leftrightarrow g$ symmetry
of the last expression.  $\Box$

\medskip
Our last lemma is the basis of the proof of part III of the theorem.
For every $V\in\mathcal{S}$ such that $|V|\geq 2$, let us define
\begin{equation}\label{P5-7}
B_V \;\triangleq\; \bigcap_{i\in V}B(i,V_i^*)\;.
\end{equation}

\begin{lemma}\label{lem3}
Let $V\in\mathcal{S}$.
\begin{fleuve11}
\item For every $W\in\mathcal{S}(\comp{V})$,
\begin{equation}\label{lem3-3} \rho_{V\cup W} \;=\;
\frac{\rho_V}{\lambda_{V}(\rho_V\, \rho_W^{-1})}
\;=\;
\frac{\rho_W}{\lambda_{W}(\rho_W\, \rho_V^{-1})}\quad \text{on } B_{V\cup W}\;.
\end{equation}
\item For every $j\in \comp{V}$
\begin{equation}\label{lem3-4}
B(j,V)\in\mathcal{F}_{\comp{V}}\;.
\end{equation}
\item If $\mu\in\mathcal{N}$, then
\begin{fleuveaa}
\item $\mu\,\lambda_V\bigl(\comp{B(k,V_k^*)}\bigr)=0$ for every $k\in V\;,$
\item $\mu\,\lambda_V(\comp{B_V})=0\;.$
\end{fleuveaa}
Furthermore, if $\mu$ satisfies the singleton consistency
\begin{equation}\label{P5} \mu \Bigl((\rho_i\, \lambda_i) (h)\Bigr)\;=\;\mu(h)
\quad\text{for every } i\in\mathbb{Z}^d\;,
\end{equation}
then
\begin{fleuveaa}
\item $\mu\bigl(\comp{B(j,V)}\bigr)=0$ for every $j\in \comp{V}\;,$
\item $\mu(\comp{B_V})=0\;.$
\end{fleuveaa}
\end{fleuve11}
\end{lemma}
\prf

{\bf (1)} Let $\omega\in B_{V\cup W}$. Then $\omega_V\in b(V,W,\omega)$
and $\omega_W\in b(W,V,\omega)$. Hence by Lemma \ref{spelem1} (3),
we have the claim.

{\bf (2)} It suffices to combine \reff{spe16}--\reff{spe19*}

{\bf (3)(a)} We apply part (2):
\begin{equation}\label{lem3-8}
\begin{aligned}
\mu\, \lambda_V\left(\comp{B(k,V_k^*)}\right)
&\;=\; \mu\, \lambda_k\,
\lambda_{V_k^*}\bigl(\comp{B(k,V_k^*)}\bigr)\nonumber\\
&\;=\; \mu\, \lambda_k\left(\comp{B(k,V_k^*)}\,
\lambda_{V_k^*}(\Omega)\right)\nonumber\\
&\;=\; 0\;,
\end{aligned}
\end{equation}
where we use the fact that the measure $\lambda^{V_k^*}$ is finite.

{\bf (3)(c)} Since $\mu\in\mathcal{N}$ satisfies SC we have
\begin{equation}\label{lem3-7}
\mu\bigl(\comp{B(j,V)}\bigr)
\;=\; \mu\,\lambda_j\bigl(\ind{\comp{B(j,V)}}\, \rho_j\bigr)
\;=\; 0\;.
\end{equation}

{\bf (3)(b)--(d)} In view of part (3)(a)--(c), the proof is a consequence of the
following
observation. For a measure $\nu$
\begin{equation}
\nu\bigl(B(k,V_k^*)\bigr)\;=\;0\;\; \forall\; k\in V
\quad \Longrightarrow \quad
\nu\bigl(\comp{B_{V}}\bigr)\;=\;0 \;.
\end{equation}
This follows from the inequality
\begin{equation}
\nu\bigl(\comp{B_V}\bigr)
\;=\; \nu\Bigl(\bigcup_{k\in V} \comp{B(k,V_k^*)}\Bigr)
\;\leq\; \sum_{k\in V}\nu\bigl( \comp{B(k,V_k^*)}\bigr).\quad \Box
\end{equation}

\paragraph{Proof of Theorem \protect \ref{th1}~:}\mbox{}
\smallskip

We consider the functions $\rho_\Lambda$ constructed in the previous
Lemma \ref{spelem1} and a bounded measurable function $h$.  We will prove, by
induction over $|\Lambda|$,
where $\Lambda \in \mathcal{S}$, that
\begin{fleuveP}
\item $\rho_{\Lambda}$ is normalized;
\item For each $\Gamma \subset \Lambda$
\begin{equation}\label{c.star}
\Bigl( \rho_{\Lambda} \, \lambda_{\Lambda}\Bigr)
\Bigl( \left( \rho_{\Gamma} \, \lambda_{\Gamma}
\right)(h)\Bigr) \;=\; (\rho_{\Lambda} \, \lambda_{\Lambda})(h)\;.
\end{equation}
\item If \reff{r.34} holds, every specification in $\Lambda$ of the
  form $\left\{\widetilde{\rho}_{\Gamma} \lambda_\Gamma : \Gamma
    \subset \Lambda \right\}$ such that
\begin{equation}\label{P3a}
\Bigl( \widetilde{\rho}_{\Lambda} \, \lambda_{\Lambda}\Bigr) \Bigl(
\left( \rho_{i} \, \lambda_{i} \right)(h)\Bigr) \;=\;
(\widetilde{\rho}_{\Lambda} \, \lambda_{\Lambda})(h)\;,\;\forall i
\in \Lambda
\end{equation}
satisfies that, for each $\omega \in \Omega_{\Lambda}$,
\begin{equation}\label{P3b}
\widetilde{\rho}_{\Lambda}(\xi_{\Lambda} \omega) =
\rho_{\Lambda}(\xi_{\Lambda} \omega) \quad \text{for }
\lambda^{\Lambda}\text{-a.a. } \xi_{\Lambda} \in \Omega_{\Lambda}\;.
\end{equation}
\item If all the functions $\rho_i$, $i\in\zzz$, are continuous and
  \reff{cont1} holds, then each function
  $\rho_\Lambda$ is continuous and for all $i \in \comp{\Lambda}$
  there exists $x_\Lambda \in \bigcap_{\omega}
  b(\Lambda,i,\omega)$ such that
\begin{equation}\label{P4}
\int \sup_\omega \left(\rho_i \, \rho_{\Lambda}^{-1}\right)
(\sigma_i x_\Lambda \omega) \, \lambda^{i}(d \sigma_i) \;<\; \infty
\;.
\end{equation}
\item If $\mu\in\mathcal{N}$ (recall \reff{N-set}) and satisfies singleton
consistency \reff{P5},
then
\begin{equation}\label{P5-2} \mu \Bigl((\rho_\Lambda\,
\lambda_\Lambda)(h)\Bigr)\;=\;\mu(h).
\end{equation}
\end{fleuveP}

The case $|\Lambda|=1$ is straightforward: (P1) is just the singleton
normalization \reff{spe33}, (P2), (P3) and (P5) are trivially true while
(P4) is \reff{cont1}.  We take now $\Lambda \in \mathcal{S}$ with
$|\Lambda| \geq 2$ and assume that (P1)--(P5) are verified by all its
non-trivial subsets.

\paragraph {Proof of (P1)}
Let $\omega \in \Omega$ and $k \in \Lambda$. By the factorization
property \reff{spe15} of $\lambda_\Lambda$ and the definition of
$\rho_\Lambda$ we have that
\begin{equation}
\lambda_{\Lambda} \left( \rho_{\Lambda} \right)(\omega) \;=\;
\lambda_{k} \left( \lambda_{\Lambda_{k}^{*}}\Bigr(
\frac{\rho_{\Lambda_{k}^{*}}}{R_{\Lambda_{k}^{*}}^{k}}\Bigr)
\right)(\omega)\;.
\end{equation}
Therefore, by the $\tribu_{\comp{(\Lambda_{k}^{*})}}$-measurability of
$R_{\Lambda_{k}^{*}}^{k}$ and the inductive normalization (P1),
\[
\lambda_{\Lambda} \left( \rho_{\Lambda} \right)(\omega) \;=\;
\lambda_{k} \left( \frac{\lambda_{\Lambda_{k}^{*}}
\left(\rho_{\Lambda_{k}^{*}}\right)}{R_{\Lambda_{k}^{*}}^{k}}\right)(\omega)
\;=\; \lambda_{k} \left(
\frac{1}{R_{\Lambda_{k}^{*}}^{k}}\right)(\omega)\;.
\]
Replacing
\[
R_{\Lambda_{k}^{*}}^{k}(\omega) \;=\;
\left(\frac{\rho_{\Lambda_{k}^{*}}}{\rho_{k}} \, \lambda_{k} \left(
    \rho_{k} \, \rho_{\Lambda_{k}^{*}}^{-1}
  \right)\right)(x_{\Lambda_{k}^{*}} \omega ) \;,
\]
for any $x_{\Lambda_{k}^{*}} \in b(\Lambda_{k}^{*},\{k\},\omega)$, we
readily obtain
$\lambda_{\Lambda} \left( \rho_{\Lambda} \right)(\omega) = 1$.
\medskip

\paragraph{Proof of (P2)~:\,\,}
It suffices to show that for some $i\in\Lambda$
\begin{equation}\label{c.10} \Bigl( \rho_{\Lambda} \, \lambda_{\Lambda} \Bigr)
  \Bigl( \left(\rho_{\Lambda_{i}^{*}} \,
    \lambda_{\Lambda_{i}^{*}}\right)(h)\Bigr) \;=\; (\rho_{\Lambda}
  \,\lambda_{\Lambda})(h)\;.
\end{equation}
Indeed, such an identity combined with the inductive hypothesis (P2)
yields that for $\Gamma$ strictly contained in $\Lambda$,
\begin{equation}\Bigl( \rho_{\Lambda} \, \lambda_{\Lambda} \Bigr) \Bigl( \left(
    \rho_{\Gamma} \, \lambda_{\Gamma} \right)(h)\Bigr) \;=\; \Bigl(
  \rho_{\Lambda} \, \lambda_{\Lambda} \Bigr) \Bigl(
  \left(\rho_{\Lambda_{i}^{*}} \,
    \lambda_{\Lambda_{i}^{*}}\right)\Bigl(\left( \rho_{\Gamma} \,
    \lambda_{\Gamma} \right)(h)\Bigr)\Bigr) \;=\; (\rho_{\Lambda} \,
  \lambda_{\Lambda})(h)\;,
\end{equation}
as needed.  To prove\reff{c.10} we use the definitions of
$\lambda_\Lambda$ and $\rho_\Lambda$ to write
\begin{equation}
\Bigl( \rho_{\Lambda} \, \lambda_{\Lambda} \Bigr) \Bigl( \left(
\rho_{\Lambda^*_i} \, \lambda_{\Lambda^*_i} \right)(h)\Bigr) \;=\;
\lambda_i \left(\lambda_{\Lambda^*_i}
\left(\frac{\rho_{\Lambda^*_i}}{R_{\Lambda^*_i}^{i}}
\lambda_{\Lambda^*_i}(\rho_{\Lambda^*_i} \, h)\right)\right)\;.
\end{equation}
Since $R_{\Lambda^*_i}^{i}$ is $\tribu_{\Lambda^{*{\rm
      c}}_i}$-measurable and
$\lambda_{\Lambda^*_i}(\rho_{\Lambda^*_i})=1$ [inductive (P1)], it
follows that
\begin{equation}\Bigl(
\rho_{\Lambda} \, \lambda_{\Lambda} \Bigr)
\Bigl( \left( \rho_{\Lambda^*_i} \, \lambda_{\Lambda^*_i}
\right)(h)\Bigr) \;=\; \lambda_i
\left(\lambda_{\Lambda^*_i}\left(\frac{\rho_{\Lambda^*_i} \, h}
{R_{\Lambda^*_i}^{i}}\right)\right)
\;=\; \left(\rho_\Lambda\,\lambda_\Lambda\right)(h)\;.
\end{equation}

\paragraph{Proof of (P3)~:\,\,}
We pick $k \in \Lambda$ and apply Lemma \ref{lem1} to the
specification 
$\left\{\widetilde{\rho}_{\Gamma} \lambda_{\Gamma}
\right\}_{\Gamma \subset \Lambda}$ 
for $f \equiv \ind{A_{\Lambda}}$ and
$g \equiv \ind{B_{\Lambda}}$ with $A_\Lambda, B_\Lambda \in
\tribu_{\Lambda}$.  We obtain
\begin{equation}
\begin{aligned}
&\int \widetilde{\rho}_{\Lambda_{k}^{*}}(\xi_{\Lambda} \omega) \,
  \widetilde{\rho}_{k}(\xi_k x_{\Lambda_{k}^{*}}
\omega) \, \widetilde{\rho}_{\Lambda}(x_\Lambda \omega) \,
  \ind{A_\Lambda}(\xi_\Lambda) \, \ind{B_\Lambda}(x_\Lambda)
\, \lambda^{\Lambda}(d \xi_{\Lambda}) \, \lambda^{\Lambda}(d
  x_{\Lambda})\hspace{-12cm}\\
& \qquad=\;
\displaystyle{\int \widetilde{\rho}_{k}(x_{\Lambda} \omega) \,
  \widetilde{\rho}_{\Lambda_{k}^{*}}(x_{\Lambda_{k}^{*}}
\xi_k \omega) \, \widetilde{\rho}_{\Lambda}(\xi_\Lambda \omega) \,
\ind{A_\Lambda}(\xi_\Lambda) \,
\ind{B_\Lambda}(x_\Lambda) \, \lambda^{\Lambda}(d x_{\Lambda}) \,
\lambda^{\Lambda}(d \xi_{\Lambda})\;,}\nonumber
\end{aligned}
\end{equation}
for every $\omega \in \Omega_{\Lambda^{\text{\rm c}}}$. Each member of
the preceding equality defines a probability measure over the product
$\sigma$-algebra $\tribu_\Lambda \otimes \tribu_\Lambda$.  This
$\sigma$-algebra is generated by the $\pi$-system, $\{A_\Lambda \times
B_\Lambda : A_\Lambda, \, B_\Lambda \in \tribu_\Lambda\}$.  As both
sides coincide on these system, they must be equal as
probability measures and, with the aid of the inductive hypothesis
(P3) we conclude that
\begin{equation}\label{uni3}
\rho_{\Lambda_{k}^{*}}(\xi_{\Lambda}
  \omega) \, \rho_{k}(\xi_k x_{\Lambda_{k}^{*}} \omega) \,
  \widetilde{\rho}_{\Lambda}(x_\Lambda \omega) \;=\;
  \rho_{k}(x_{\Lambda} \omega) \,
  \rho_{\Lambda_{k}^{*}}(x_{\Lambda_{k}^{*}} \xi_k \omega) \,
  \widetilde{\rho}_{\Lambda}(\xi_\Lambda \omega)\;,
\end{equation}
for $\lambda^\Lambda \times \lambda^\Lambda$-a.a.\ $(\xi_\Lambda,
x_\Lambda) \in \Omega_\Lambda \times \Omega_\Lambda$.  Since by
assumption each $\lambda^j$ charges $b(j,\Lambda_{j}^{*},\omega)$,
identity \reff{uni3} must be verified for some choice of $x_j \in
b(j,\Lambda_{j}^{*},\omega)$.  In this case the factors of
$\widetilde\rho_\Lambda$ in the RHS of \reff{uni3} are non-zero and we
can solve
\begin{equation}
  \widetilde{\rho}_{\Lambda}(\xi_\Lambda \omega) \;=\;
  \frac{\rho_{\Lambda_{k}^{*}}(\xi_{\Lambda} \omega) \, \rho_{k}(\xi_k
    x_{\Lambda_{k}^{*}} \omega) \,
    \widetilde{\rho}_{\Lambda}(x_\Lambda \omega)}{\rho_{k}(x_{\Lambda}
    \omega) \, \rho_{\Lambda_{k}^{*}}(x_{\Lambda_{k}^{*}} \xi_k
    \omega)}
\end{equation}
for $\lambda^\Lambda$-a.a.\ $\xi_\Lambda \in \Omega_\Lambda$.  If we
integrate both sides with respect $\lambda^\Lambda(\xi_\Lambda)$, we
get
\begin{equation}
1 \;=\;
\Bigr(\lambda_{\Lambda}(\widetilde{\rho}_{\Lambda})\Bigr)(\omega)
  \;=\; \frac{\widetilde{\rho}_{\Lambda}(x_\Lambda
    \omega)}{\rho_{k}(x_{\Lambda} \omega)}\, \lambda_k \left( \rho_{k}
    \, \rho_{\Lambda_{k}^{*}}^{-1} \,
    \lambda_{\Lambda_{k}^{*}}(\rho_{\Lambda_{k}^{*}})\right)(x_{\Lambda}
  \omega)\;.
\end{equation}
Since
$\lambda_{\Lambda_{k}^{*}}(\rho_{\Lambda_{k}^{*}})\equiv 1$, we obtain
\begin{equation}\label{uni12}
\widetilde{\rho}_{\Lambda}(x_\Lambda \omega) \;=\;
  \frac{\rho_{k}(x_{\Lambda} \omega)}{\lambda_k \left( \rho_{k} \,
      \rho_{\Lambda_{k}^{*}}^{-1}\right)(x_{\Lambda} \omega)}
\;=\; \rho_{\Lambda}(x_\Lambda \omega)\;.
\end{equation}
From \reff{uni3} and \reff{uni12}, we conclude that each
$\widetilde{\rho}_\Lambda$ satisfying \reff{P3a} is
$\lambda^{\Lambda}$-a.s.\ uniquely determined. Since $\rho_\Lambda$
itself satisfies \reff{P3a}, statement \reff{P3b} follows.

\paragraph{Proof of (P4)~:\,\,}
We first remark that if $V\subset\comp{\Lambda}$ we can construct some
$x_\Lambda\in\bigcap_\omega b(\Lambda,V,\omega)$ simply by choosing
$x_j\in\bigcap_\omega b(j,V\cup\Lambda^*_j,\omega)$ [see definition
\reff{spe13}].  Let $k \in \Lambda$ and $x_{\Lambda_{k}^{*}} \in
\bigcap_{\omega}b(\Lambda_{k}^{*},k,\omega)$. The inductive hypotheses
(P4) implies the continuity of the functions $(\rho_k \,
\rho_{\Lambda_{k}^{*}}^{-1}) (\sigma_k x_{\Lambda_{k}^{*}} \,\cdot\,)$
for each $\sigma_k\in E$.  These functions are uniformly bounded above
by $\sup_{\omega} \rho_k \, \rho_{\Lambda_{k}^{*}}^{-1} (\sigma_k
x_{\Lambda_{k}^{*}} \omega)$ which ---by the inductive assumption
\reff{P4}--- is integrable with respect to $\lambda^{k}(d\sigma_k)$.
The sequential continuity of the function $\lambda_k (\rho_k \,
\rho_{\Lambda_{k}^{*}}^{-1})(x_{\Lambda_{k}^{*}}\, \cdot \,)$ follows,
then, from the dominated convergence theorem.  This function is
strictly positive because of the choice of $x_{\Lambda_{k}^{*}}$.
These continuity and non-nullness, plus the inductive continuity
hypothesis, imply that
\begin{equation}
\rho_{\Lambda}(\,\cdot\,) \;\triangleq\;
    \frac{\rho_{\Lambda_{k}^{*}}(\,\cdot\,) \,
      \rho_k(x_{\Lambda_{k}^{*}}
      \,\cdot\,)}{\rho_{\Lambda_{k}^{*}}(x_{\Lambda_{k}^{*}}
      \,\cdot\,) \, \lambda_k (\rho_k \,
      \rho_{\Lambda_{k}^{*}}^{-1})(x_{\Lambda_{k}^{*}} \,\cdot\,)}
\end{equation}
is a continuous function.

Finally we prove \reff{P4}.  The existence of some $x_\Lambda \in
\bigcap_{\omega} b(\Lambda, i, \omega)$ yields the identity
\begin{equation}
\left(\rho_i \, \rho_{\Lambda}^{-1}\right)
  (\sigma_i x_\Lambda \omega) \;=\;
\left(\rho_i \, \rho_{k}^{-1}\right)(\sigma_i
  x_\Lambda \omega) \times \int \left(\rho_k \,
  \rho_{\Lambda_{k}^{*}}^{-1}\right)
(\sigma_k \sigma_i x_\Lambda \omega) \,
  \lambda^{k}(d \sigma_k)\; ,
\end{equation}
valid for \emph{all} $\omega \in \Omega_{\comp{(\Lambda \cup
    \{i\})}}$, each $k \in \Lambda$ and each $\sigma_i \in \Omega_i$.
We take supremum over $\omega$ and integrate with respect to
    $\lambda^{i}$ to obtain
\begin{equation}
\begin{aligned}
&\int \sup_\omega \left(\rho_i \, \rho_{\Lambda}^{-1}\right)
 (\sigma_i x_\Lambda \omega)\,  \lambda^{i}(d \sigma_i)\\
 &\qquad \displaystyle{\;\leq\; \int \sup_\omega \Bigl(\rho_i \,
    \rho_{k}^{-1}\Bigr)(\sigma_i x_\Lambda \omega) \, \lambda^{i}(d
    \sigma_i) \times \int \sup_\omega \left( \rho_k \,
\rho_{\Lambda_{k}^{*}}^{-1}\right)(\sigma_k
    x_{\Lambda^*_k} \omega) \, \lambda^{k}(d \sigma_k)\;.}\nonumber\\
\end{aligned}
\end{equation}
Both integrals in the RHS are finite by the inductive assumption (P4).

\paragraph {Proof of (P5)~:\,\,}
Fix $i\in\Lambda$. Since $\mu\in\mathcal{N}$ satisfies singleton consistency
\reff{P5},
\begin{equation}\label{P5-13}
\begin{aligned}
\mu\bigl((\rho_\Lambda\, \lambda_\Lambda)(h)\bigr)
&\;=\; \mu\left(\lambda_\Lambda\left(\ind{B_\Lambda}\,
\rho_\Lambda\, h\right)\right)\nonumber\\
&\;=\; \mu\left(\lambda_\Lambda\left(\ind{B_\Lambda}
\frac{\rho_i}{\lambda_i\bigl(\rho_i\,\rho_{\Lambda_i^*}^{-1}\bigr)}\,
h\right)\right)\;,
\end{aligned}
\end{equation}
where the first and second identities come respectively from parts Lemma (3)(b)
and (1) of Lemma \ref{lem3}. We write $B_\Lambda =
\bigcap_{k\in\Lambda}B(k,\Lambda_k^*)$
and decompose $\lambda_\Lambda=\lambda_{\Lambda_i^*}\, \lambda_i$. Using the
measurability
and the support property of parts (2) and (3)(b) of Lemma \ref{lem3}, we see
that
\begin{equation}
\begin{aligned}
\label{P5-15} \mu\bigl((\rho_\Lambda\, \lambda_\Lambda)(h)\bigr)
&\;=\; \mu\left(\lambda_{\Lambda_i^*}\left(\ind{}\Bigl\{
\bigcap_{j\in \Lambda_i^*}B(j,\Lambda_j^*)\Bigr\}
\frac{\displaystyle{\lambda_i\bigl(\ind{B(i,\Lambda_i^*)}\, \rho_i\, h\bigr)}}
{\lambda_i\bigl(\rho_i\,\rho_{\Lambda_i^*}^{-1}\bigr)}\right)\right)\nonumber\\
&\;=\; \mu\left(\ind{B(i,\Lambda_i^*)}\lambda_{\Lambda_i^*}\left(
\ind{}\Bigl\{\bigcap_{j\in \Lambda_i^*}B(j,\Lambda_j^*)\Bigr\}
\frac{\displaystyle{\lambda_i\bigl(\ind{B(i,\Lambda_i^*)}\, \rho_i\, h\bigr)}}
{\lambda_i\bigl(\rho_i\,\rho_{\Lambda_i^*}^{-1}\bigr)}\right)\right)\nonumber\\
&\;=\; \mu\left(\lambda_{\Lambda_i^*}\left(\ind{B_\Lambda}
\frac{\displaystyle{\lambda_i\bigl(\ind{B(i,\Lambda_i^*)}\, \rho_i\, h\bigr)}}
{\lambda_i\bigl(\rho_i\,\rho_{\Lambda_i^*}^{-1}\bigr)}\right)\right)\;.
\end{aligned}
\end{equation}
The inductive hypotheses (P5) implies that $\mu(\rho_{\Lambda_i^*}\,
\lambda_{\Lambda_i^*})=\mu$. Hence
\begin{equation}
\begin{aligned}
\label{P5-19} \mu\bigl((\rho_\Lambda\, \lambda_\Lambda)(h)\bigr)
&\;=\; \mu\left(\frac{\ind{B_\Lambda}}{\rho_{\Lambda_i^*}}
\frac{\displaystyle{\lambda_i\bigl(\ind{B(i,\Lambda_i^*)}\, \rho_i\, h\bigr)}}
{\lambda_i\bigl(\rho_{i}\,\rho_{\Lambda_i^*}^{-1}\bigr)}\right)\nonumber\\
&\;=\; \mu\left(\frac{1}{\rho_{\Lambda_i^*}}
\frac{\displaystyle{\lambda_i\bigl(\ind{B(i,\Lambda_i^*)}\, \rho_i\, h\bigr)}}
{\lambda_i\bigl(\rho_{i}\,\rho_{\Lambda_i^*}^{-1}\bigr)}\right)\;,
\end{aligned}
\end{equation}
where the second line comes from support property of Lemma \ref{lem3} (3)(d).
By singleton consistency \reff{P5}
\begin{equation}
\begin{aligned}
\label{P5-20}
\mu\bigl((\rho_\Lambda\, \lambda_\Lambda)(h)\bigr)
&\;=\; \mu\,\lambda_i\left(\frac{\rho_i}{\rho_{\Lambda_i^*}}
\frac{\displaystyle{\lambda_i\bigl(\ind{B(i,\Lambda_i^*)}\, \rho_i\, h\bigr)}}
{\lambda_i\bigl(\rho_{i}\,\rho_{\Lambda_i^*}^{-1}\bigr)}\right)\nonumber\\
&\;=\; \mu\bigl(\lambda_i(\ind{B(i,\Lambda_i^*)}\, \rho_i\,
h)\bigr)\;.
\end{aligned}
\end{equation}
Therefore
\begin{equation}
\label{P5-21} \mu\bigl((\rho_\Lambda\, \lambda_\Lambda)(h)\bigr)
\;=\; \mu\bigl(\lambda_i(\rho_i\, h)\bigr)
\;=\; \mu(h)\;,
\end{equation}
where once again we use Lemma \ref{lem3} (3)(b) and singleton consistency.
$\quad \Box$

\section{Proof of Proposition \ref{prop1}}
The proof relies on results already stated in Georgii (1988).
Since $\lambda^i(b(i,\{j\},\alpha))$ is strictly positive
for all $\alpha\in\Omega$
and $i\neq j \in\mathbb{Z}^d$, we can apply Proposition (1.30) of
Georgii (1988) to conclude that, for
$\lambda_{\{i,j\}}(\,\cdot\mid\alpha)$-almost all
$\omega\in\Omega$,
\begin{equation}\label{pr-03} \rho_{\{i,j\}}(\omega)\,
\rho_i(x_i\omega) \;=\; \rho_{\{i,j\}}(x_i\omega)\, \rho_i(\omega)
\end{equation}
and
\begin{equation}\label{pr-05} \rho_{\{i,j\}}(\omega)\,
\rho_j(x_j\omega) \;=\; \rho_{\{i,j\}}(x_j\omega)\, \rho_j(\omega)
\end{equation}
for all $x_i\in b(i,\{j\},\omega)$ and $x_j\in b(j,\{i\},\omega)$.
These identities imply, by the definition (\ref{spe19}) of good
sets, that for $\lambda_{\{i,j\}}(\,\cdot\mid\alpha)$-almost all
$\omega\in\Omega$,
\begin{equation} \label{pr-06}
\frac{\rho_{\{i,j\}}(x_j\omega)\,
\rho_j(\omega)}{\rho_j(x_j\omega)}
\;=\;\frac{\rho_{\{i,j\}}(x_i\omega)\,
\rho_i(\omega)}{\rho_i(x_i\omega)}
\end{equation}
for all $x_i\in
b(i,\{j\},\omega)$ and $x_j\in b(j,\{i\},\omega)$. But, by Theorem
(1.33) of Georgii (1988), we have that for these $\omega$, $x_i$
and $x_j$,
\begin{equation}\label{pr-07}
\rho_{\{i,j\}}(x_j\omega)
\;=\;\frac{\rho_i(x_j\omega)}{\lambda_i(\rho_i\,\rho_j^{-1})(x_j\omega)}
\end{equation}
and
\begin{equation}\label{pr-09}
\rho_{\{i,j\}}(x_i\omega)
\;=\;\frac{\rho_j(x_i\omega)}{\lambda_j(\rho_j\,\rho_i^{-1})(x_i\omega)}.
\end{equation}
The substitution of (\ref{pr-07}) and (\ref{pr-09}) into
(\ref{pr-06}) yields the result. In particular when $E$ is
countable and each $\lambda^i$ is the counting measure, as a
consequence of Proposition (1.30) of Georgii (1988),
(\ref{pr-03}--\ref{pr-05}) hold for all $\omega\in\Omega$, $x_i\in
b(i,\{j\},\omega)$ and $x_j\in b(j,\{i\},\omega)$. Thus in that
case, (H2) is fulfilled for all $\omega\in\Omega$. $\quad
\Box$

\section*{Acknowledgements}  G.M. gratefully acknowledges support from
the Netherlands Organization for Scientific Research during his stay at
EURANDOM.
The authors would like to thank the anonymous referee for his careful comments.

\bibliographystyle{alea2}

\begin{thebibliography}{}

\bibitem{dacnah01}
S.~Dachian and B.~S.\ Nahapetian.
\newblock Description of random fields by means of one-point conditional
  distributions and some applications.
\newblock {\em Markov Proc.\ Rel.\ Fields}, 7:193--214, 2001.

\bibitem{dacnah04}
S.~Dachian and B.~S.\ Nahapetian.
\newblock Description of specifications by means of probability distributions
  in small volumes under condition of very weak positivity.
\newblock {\em J. Stat. Phys.}, 117:281--300, 2004.

\bibitem{dob68}
R.~L.\ Dobrushin.
\newblock The description of a random field by means of conditional
  probabilities and conditions of its regularity.
\newblock {\em Theory of probability and its applications}, 13:197--224, 1968.

\bibitem{vEFS_JSP}
A.~C.~D.\ van Enter, R.~Fern{\'a}ndez and A.~D.\ Sokal.
\newblock Regularity properties and pathologies of position-space
  renormalization-group transformations: scope and limitations of {G}ibbsian
  theory.
\newblock {\em J. Stat. Phys.}, 72:879--1167, 1993.

\bibitem{entmaeshl00}
\newblock A.~C.~D.\ van~Enter, C.~Maes and S.~Shlosman.
\newblock Dobrushin's program on {G}ibbsianity restoration:  {W}eakly
{G}ibbs and almost {G}ibbs random fields.
\newblock {\em Amer.\ Math.\ Soc.\ Transl.}, 198:59--69, 2000.

\bibitem{fermai03a}
 R.~Fern\'andez and G.~Maillard.
 \newblock Chains with complete connections and one-dimensional {G}ibbs
   measures.
 \newblock {\em Electron.\ J.\ Probab.}, 9:145--76, 2004.

 \bibitem{flosul80}
 R.~G. Flood and W.~G. Sullivan.
 \newblock Consistency of random fields specifications.
 \newblock {\em Z. Wahrsch. Verw. Gebiete}, 53(2):147--156, 1980.

\bibitem{geo88}
H.-O.\ Georgii.
\newblock {\em Gibbs Measures and Phase Transitions}.
\newblock Walter de Gruyter (de Gruyter Studies in Mathematics, Vol.\ 9),
  Berlin--New York, 1988.

\bibitem{koz74}
O.~K.\ Kozlov.
\newblock {G}ibbs description of a system of random variables.
\newblock {\em Probl. Inform. Transmission}, 10:258--65, 1974.

\bibitem{pres04}
C.\ Preston.
\newblock {O}ne point Gibbs states.
\newblock {\em Preprint}, May 2004.
\newblock It used to be available at \texttt{ http://www.mathematik.uni-bielefeld.de/$\sim$preston/}

\bibitem{sok81}
A.~D.\ Sokal.
\newblock Existence of compatible families of proper regular
conditional probabilities.
\newblock {\em Z.\ Wahrscheinlichkeitstheorie verw.\ Gebiete},
56:537--48, 1981.

\bibitem{sul73}
W.~G.\ Sullivan.
\newblock {P}otentials for almost {M}arkovian random fields.
\newblock {\em Comm. Math. Phys.}, 33:61--74, 1973.


\end{thebibliography}

\end{document}